%% file: Niu-Wang-Xie-2.tex
\documentclass[review,hidelinks,onefignum,onetabnum]{siamart220329}

\input{ex_shared}

\ifpdf
\hypersetup{
  pdftitle={Global dynamics of three-dimensional Lotka-Volterra competition models with seasonal succession: II. Uniqueness of positive fixed points},
  pdfauthor={Lei Niu, Yi Wang, and Xizhuang Xie}
}
\fi

\begin{document}

\maketitle

\begin{abstract}
In this second part of the series, we investigate the uniqueness of positive fixed points of the Poincar\'e map associated with the 3-dimensional Lotka-Volterra competition model with seasonal succession. Building on our first part of the series on the classification of 33 dynamical equivalence classes (regardless of the uniqueness of  positive fixed points), we demonstrate in this paper that classes 26 and 27 may indeed exhibit multiple positive fixed points. This reveals a fundamental distinction from both its 2-dimensional analogue and the classical 3-dimensional competitive Lotka-Volterra model. More concretely, by focusing on the model with identical growth and death rates, we establish an equivalent characterization for the (non)uniqueness of positive fixed points. Based on this characterization, we further show that classes 19–25 and 28–33 admit a unique positive fixed point and exhibit trivial dynamics: all trajectories converge to some fixed point (corresponding to harmonic solutions). In contrast, classes 26 and 27 possess richer dynamical scenarios:  there can contain a continuum of invariant closed curves, on which orbits may be periodic (corresponding to subharmonic solutions), dense (corresponding to quasi-periodic solutions), or may even consist entirely of positive fixed points (which exhibits the nonuniqueness of positive fixed points).

\end{abstract}



\section{Statement of results}\label{section:1}
\subsection{Introduction}
This is the second of the series. In Part I \cite{nwx2023a}, we investigated the long-term
behaviour of the three-dimensional
Lotka-Volterra competition model with seasonal
succession
\begin{equation}\label{seasonal-system}
 \begin{dcases}
\frac{d x_{i}}{d t}=-\mu_{i} x_{i}, \quad t\in[k \omega,  k \omega+(1-\varphi) \omega), \\
\frac{d x_{i}}{d t}=x_{i}\Big(b_{i}-\sum_{j=1}^{3} a_{i j} x_{j}\Big),\quad t\in [k \omega+(1-\varphi) \omega,  (k+1) \omega),
\end{dcases}
\end{equation}
where $k \in \mathbb{Z}_{+}, \varphi \in (0,1]$ and $\omega$, $\mu_{i}$, $b_{i}$, and $a_{ij}$ are all positive constants, for  $i,j=1,2,3$.  The parameter $\omega$ is the length of a full seasonal cycle, and $\varphi$ is the fraction of that cycle occupied by the “good” season.

Following Klausmeier \cite{k2010,k2012}, system \eqref{seasonal-system} serves as a prototype for studying how seasonal alternation shapes competitive outcomes. It is known that if the average growth rate $$
r_i:=b_{i} \varphi -\mu_{i}(1-\varphi)\leq 0,
$$
the species $i$ is driven to extinction. Conversely, when \(r_i > 0\) for all species, the associated Poincar\'e map $\mathcal{P}$ of \eqref{seasonal-system} possesses a global attractor $\Sigma_\mathcal{P}\subset\mathbb{R}^3_+:=\{x\in \mathbb{R}^3:x_i\geq 0,~i=1,2,3\}$, known as the carrying simplex. This invariant, codimension-one manifold attracts every nontrivial trajectory \cite{baigent2019,diekmann2008carrying, hirsch1988, hou2020, LG, ruiz2011exclusion,Wang-Jiang-02}.

However, characterizing the dynamics on $\Sigma_\mathcal{P}$ remains a significant challenge.
To the best of our knowledge, almost all previous work in this field has relied on the assumption of uniqueness of positive fixed points; see, for example, \cite{jiang-niu2017,jnz,Z993} for autonomous competitive ODEs and \cite{Gyllenberg2019,Gyllenberg2020b,jiang2014,LG} for competitive maps.
Hsu and Zhao \cite{Hsu-Zhao} proved that such uniqueness holds automatically for the two-species version of model \eqref{seasonal-system} using the theory of monotone dynamical systems. They further
classified the global dynamics of the two-species model in \cite{Hsu-Zhao}, which closely resemble the
classical dynamic scenarios in the Lotka-Volterra competition model.
In recent studies, more two-species competition models that account for seasonal succession have been proposed and examined.  Examples include the competition model between Aedes aegypti and Aedes albopictus introduced in \cite{zheng2023,Zheng2025}, as well as models analyzing competitive dynamics under harvesting strategies in \cite{yu2024}.


Yet, for the three-species version of model \eqref{seasonal-system}, the uniqueness of the positive fixed point of
\(\mathcal{P}\) remains a complex and poorly understood problem. The main difficulty stems from the fact that, as a typical example of time-periodically forced differential equations, model \eqref{seasonal-system} lacks an explicit expression for  \(\mathcal{P}\).

In Part I \cite{nwx2023a}, without assuming the uniqueness of positive fixed points, we introduced an equivalence relation for all Poincaré maps associated with model \eqref{seasonal-system}, relative to the boundary dynamics. By developing an index formula on \(\Sigma_\mathcal{P}\) that reveals the intrinsic relationship among fixed points, we classified the dynamics of the model into 33 classes. Among these, classes 1–18 have no positive fixed point, with all orbits converging to a boundary fixed point. In contrast, classes 19–33 have at least one positive fixed point, although the uniqueness of these points is not guaranteed. Notably, in classes 19–25, every orbit converges to some fixed point, while the positive fixed point in class 33 is globally asymptotically stable, assuming that \(\mathcal{P}\) has a unique positive fixed point in these classes.

A natural question then arises: \emph{Could the mapping $\mathcal{P}$ admit non-unique positive fixed points?} If so, the system would exhibit a fundamental difference, distinguishing it from both its 2-dimensional analogue and the classical 3-dimensional competitive Lotka-Volterra model (see Zeeman \cite{Z993}). This is the focus of the current part of the series. Indeed, it is surprisingly difficult to construct explicit examples with multiple positive fixed points, even in the 3-dimensional case. We therefore consider in this paper a special yet general case in which all species share the same intrinsic growth rate during the good season and the same death rate during the bad season given by
\begin{equation}\label{seasonal-system-equal-db}
\begin{dcases}
\frac{d x_{i}}{d t} = -\mu x_{i}, \quad t \in [k \omega, k \omega + (1 - \varphi) \omega), \\
\frac{d x_{i}}{d t} = x_{i} \left(b - \sum_{j=1}^{3} a_{i j} x_{j}\right), \quad t \in [k \omega + (1 - \varphi) \omega, (k+1) \omega),
\end{dcases}
\end{equation}
where \(b_i = b > 0\) and \(\mu_i = \mu > 0\) for all \(i = 1, 2, 3\), and the average growth rate
\[
r := b \varphi - \mu(1 - \varphi) > 0.
\]

Throughout the paper we assume the interaction matrix $A=(a_{ij})_{3\times 3}$ is nonsingular, ensuring that the corresponding 3-dimensional autonomous Lotka-Volterra competition
model
\begin{equation}\label{PEVL}
\frac{d x_i}{dt}=x_i\Big(b-\sum_{j=1}^3 a_{ij}x_j\Big),\qquad i=1,2,3,
\end{equation}
admits at most one positive equilibrium. Under these assumptions we derive necessary and sufficient conditions for both the uniqueness and the nonuniqueness of positive fixed points of the Poincar\'e map associated with \eqref{seasonal-system-equal-db}.

\subsection{Precise statement of results}
Denote by \(\Phi_t(x)\) the solution flow of system \eqref{PEVL} with initial state \(x \in \mathbb{R}_+^3\). We now present the main results of this paper, which are proven in Section \ref{sec:3}.

\begin{theorem}[(Non)uniqueness of positive fixed points]\label{thm:unique-PFP}
For system \eqref{seasonal-system-equal-db} with seasonal succession, the following conclusions hold:
\begin{itemize}
\item[{\rm (a)}] \(\mathcal{P}\) has no positive fixed point if and only if \(\Phi_t\) of system \eqref{PEVL} has no positive equilibrium.
\item[{\rm (b)}] \(\mathcal{P}\) has a unique positive fixed point if and only if \(\Phi_t\) has positive equilibrium and satisfies:
\begin{itemize}
\item[{\rm (i)}] \(\Phi_t\) has no positive periodic orbit; or
\item[{\rm (ii)}] any positive periodic orbit \(\Gamma\) {\rm(}with minimal positive period \(T_\Gamma\){\rm)} of \(\Phi_t\) satisfies that
\[
\eta \triangleq \left(\frac{\omega}{T_\Gamma} : \frac{b}{r}\right) \notin \mathbb{Z}.
\]
\end{itemize}
\item[{\rm (c)}] \(\mathcal{P}\) has multiple positive fixed points if and only if \(\Phi_t\) has a periodic orbit \(\Gamma\) with \(\eta \in \mathbb{Z}\).
\end{itemize}
\end{theorem}

Theorem \ref{thm:unique-PFP} establishes an equivalent characterization for the (non)uniqueness of positive fixed points, by connecting with the positive periodic orbits of the flow $\Phi_t$ for the autonomous system (\ref{PEVL}). Based on such characterization, we further delineate all possible global dynamical scenarios of the system in the following theorems.

\begin{theorem}[Equivalence classes with unique positive fixed point]\label{thm_class_1}\label{thm:1.2}
For system \eqref{seasonal-system-equal-db} with seasonal succession,
among the $33$ equivalence classes {\rm (}established in {\rm \cite{nwx2023a}}{\rm )}, one has
\begin{enumerate}[\rm (i)]
    \item Classes $1$--$18$: each class admits no positive fixed point, and every orbit converges to a fixed point.
    \item Classes $19$--$25~\&~28$--$33$: each class possesses a unique positive fixed point, and every orbit converges to a fixed point.

\end{enumerate}
\end{theorem}

Theorem \ref{thm:1.2}, in a detailed version, will be proved in Section 3 (see Theorems \ref{thm:whole-dynamics} and \ref{thm:class_28_33}, respectively). We emphasize that for classes 28-32, the global convergence dynamics are studied
by exploiting the delicate interplay between the associated Poincar\'e map of system \eqref{seasonal-system-equal-db} and the autonomous flow of \eqref{PEVL}. It should be further noted that the global dynamics of classes 28–32 remain open for the discrete-time competitive mappings studied in \cite{Gyllenberg2020a, Gyllenberg2020b, jiang2014, LG, MNR2019, Niu-Ruiz-2018}, although the uniqueness of positive fixed points is already established in such setting.


\begin{theorem}[Equivalence classes admitting multiple positive fixed points]\label{thm_class_2}\label{thm:1.3} For system \eqref{seasonal-system-equal-db} with seasonal succession, classes $26$ and $27$ {\rm (}established in {\rm \cite{nwx2023a}}{\rm )} can each exhibit a continuum of invariant closed curves in ${\rm Int}\mathbb{R}^3_+$. On such a curve, one of the following alternatives exactly occurs: {\rm(}a{\rm)} any orbit is a fixed point; {\rm(}b{\rm)} any orbit is periodic; {\rm(}c{\rm)} any orbit is dense in such curve.

\end{theorem}

Theorem \ref{thm:1.3}, in a detailed version, will be proved in Section 3 (see Theorems \ref{prop:class_26} and \ref{prop:class_27}, respectively).
Together with Theorem \ref{thm:1.2},
these results establish that for system \eqref{seasonal-system-equal-db}, classes 19–25 and 28–33 possess a unique positive fixed point, whereas classes 26 and 27 may indeed admit multiple positive fixed points. This outcome highlights a fundamental distinction between the 3-dimensional system \eqref{seasonal-system} and both its 2-dimensional counterpart studied in \cite{Hsu-Zhao} and classical 3-dimensional competitive maps \cite{Gyllenberg2020a, Gyllenberg2020b, jiang2014, LG, ruiz2011exclusion}, the latter of which can only admit at most one positive fixed point. Moreover, in classes 26 and 27 for system \eqref{seasonal-system-equal-db}, when invariant closed curves appear, the corresponding positive fixed points, periodic orbits, or dense orbits on these curves yield, respectively, infinitely many harmonic, subharmonic, or quasi-periodic positive solutions to the system.

\section{Preliminaries}\label{section:2}
We denote by $\dot{\mathbb{R}}^3_+:= \{x\in \mathbb{R}^3_+: x_i>0, i=1,2,3\}$ the interior of $\mathbb{R}^3_+$ and $\partial \mathbb{R}^3_+:=\mathbb{R}^3_+\setminus \dot{\mathbb{R}}^3_+$ the boundary of $\mathbb{R}^3_+$. For each nonempty $I\subset \{1,2,3\}$, we set $H_I^+=\{x\in \mathbb{R}_+^3: x_j=0\ \mathrm{for}\  j \notin I\}$ and $\dot{H}_I^+=\{x\in H_I^+: x_i>0 \ \mathrm{for}\  i \in I\}$. In particular, $H^+_{\{i\}}$ denotes the $i$th positive coordinate axis for $i\in I$. The symbol $0$  stands for both the origin of $\mathbb{R}^3$ and the real number $0$.

The system \eqref{seasonal-system} admits a unique solution $\Psi(t, x)$ through any initial point $x\in \mathbb{R}^3$ for all $\varphi \in (0,1]$ and positive $\omega$, $\mu_{i}$, $b_{i}$ and $a_{ij}$. Clearly, the domain of $\Psi(\cdot, x)$ includes $[0,+\infty)$
in case $x\in \mathbb{R}_+^3$.
Since the system is $\omega$-periodic, we consider the associated Poincar\'{e} map $\mathcal{P}$ given by
$$\mathcal{P}(x)=\Psi(\omega, x)$$
for those $x$ for which the right-hand side is defined, namely an open set $W\subset \mathbb{R}^3$ containing $\mathbb{R}_+^3$. For each nonempty $I\subset \{1,2,3\}$, $H_I^+$ and $\dot{H}_I^+$ are positively invariant under $\mathcal{P}$.

\begin{lemma}[Lemma 2.1 in \cite{Hsu-Zhao}]
If $r_i\leq 0$, then $(\mathcal{P}_i|_{H_{\{i\}}^+})^k(x)\to 0$ as $k\to +\infty$ for all $x\in H_{\{i\}}^+$, while there is a unique nontrivial fixed point $q_{\{i\}}\in  H_{\{i\}}^+$ which is globally attracting in $H_{\{i\}}^+\setminus \{0\}$ if $r_i> 0$.
\end{lemma}

\begin{lemma}[\cite{Niu-Wang-Xie}]\label{carrying-simplex}
If $r_i>0$, $i=1,2,3$, then the Poincar\'{e} map $\mathcal{P}$ admits a carrying simplex $\Sigma_\mathcal{P}$, that is a subset  of $\mathbb{R}_+^3\setminus \{0\}$ with the  properties:
\begin{enumerate}[{\rm(P1)}]
\item $\Sigma_\mathcal{P}$ is compact, invariant and unordered;
\item $\Sigma_\mathcal{P}$ is homeomorphic via radial projection to the $2$-dimensional standard probability simplex $\Delta^{2}=\{x\in \mathbb{R}_+^3:x_1+x_2+x_3=1\}$;
\item for any $x\in \mathbb{R}_+^3\setminus \{0\}$, there exists some $z\in \Sigma_\mathcal{P}$ such that
$$\lim_{k\to \infty}\|\mathcal{P}^k(x)-\mathcal{P}^k(z)\|=0.$$
\end{enumerate}
\end{lemma}
We denote the boundary of the carrying simplex $\Sigma_\mathcal{P}$ relative to $\mathbb{R}^3_+$ by $\partial \Sigma_\mathcal{P}=\Sigma_\mathcal{P}\cap \partial\mathbb{R}^3_+$  and the interior of $\Sigma_\mathcal{P}$ relative to $\mathbb{R}^3_+$ by $\dot{\Sigma}_\mathcal{P}=\Sigma_\mathcal{P} \setminus \partial \Sigma_\mathcal{P}$.

\begin{lemma}[\cite{nwx2023a}]\label{lemma-PF}
A point $\theta \in \dot{\mathbb{R}}_{+}^{3}$ is a positive fixed point of $\mathcal{P}$ if and only if $\hat{\theta}:=\int_{0}^{\varphi \omega} \Phi_{t}(L \theta) d t$ is a positive solution of the linear algebraic system $(Ax^\tau)_i=r_i$ with $r_i>0$, $i=1,2,3$, where $A$ is the ${3\times 3}$ matrix with entries $a_{ij}$ given in the system \eqref{seasonal-system}.
\end{lemma}

\begin{lemma}[\cite{nwx2023a}]\label{prop-lambda}
For any $x\in \dot{ \mathbb{R}}^3_+$, 	there exists $(D\mathcal{P}(x))^{-1}$ which is a positive matrix.
Moreover, for any $\theta \in \mathrm{Fix}(\mathcal{P},\dot{ \mathbb{R}}^3_+)$, the eigenvalue of $D\mathcal{P}(\theta)$ with the smallest modulus, say $\lambda$, satisfies $0<\lambda<1$.
\end{lemma}

Next we recall the classification provided in \cite{nwx2023a} for all the associated Poincar\'{e} maps of the 3-dimensional systems \eqref{seasonal-system} with a carrying simplex, which is denoted by
$$
\mathrm{CLVS}(3):=\left\{\begin{array}{l}
    \mathcal{P}\in \mathcal{T}(\mathbb{R}_+^3): \mathcal{P}(x)=\Psi(\omega, x),~\Psi{~ \text{is~the solution~map~of}~}\eqref{seasonal-system}  \\
    \noalign{\smallskip}
      \text{with~} 0<\varphi<1,\omega,\mu_i,b_i,r_i,a_{ij}>0, i,j=1,2,3
\end{array}
\right\}.
$$

On the boundary of its carrying simplex, the Poincar\'e map $\mathcal{P}\in \mathrm{CLVS}(3)$ has three axial fixed points, say
$q_{\{1\}}=(q_1,0,0)$, $q_{\{2\}}=(0,q_2,0)$, $q_{\{3\}}=(0,0,q_3)$,
besides the trivial fixed point $0$, and it may have a planar fixed point $v_{\{k\}}$ in the interior of the coordinate plane $\Pi_k$. Two maps $\mathcal{P}, \hat{\mathcal{P}} \in \mathrm{CLVS}(3)$ are said to be {\it equivalent relative to the boundary of the carrying simplex} if there exists a permutation $\sigma$ of $\{1,2,3\}$ such that $\mathcal{P}$ has
a fixed point $q_{\{i\}}$ {\rm(}or $v_{\{k\}}${\rm)} if and only if $\hat{\mathcal{P}}$ has a fixed
point $\hat{q}_{\{\sigma(i)\}}$ {\rm(}or $\hat{v}_{\{\sigma(k)\}}${\rm)}, and
further $q_{\{i\}}$ {\rm(}or $v_{\{k\}}${\rm)} has the same hyperbolicity and local dynamics on the carrying simplex as $\hat{q}_{\{\sigma(i)\}}$ {\rm(}or $\hat{v}_{\{\sigma(k)\}}${\rm)}. $\mathcal{P}\in\mathrm{CLVS}(3)$ is said to be {\it stable relative to the boundary of the carrying simplex} if
all the fixed points on $\partial \Sigma_\mathcal{P}$ are hyperbolic. An
equivalence class is {\it stable} if each mapping in
it is stable relative to the boundary of the carrying simplex.

\begin{lemma}\label{lemma:PFP}
   There are a total of $33$ stable equivalence classes for $\mathrm{CLVS}(3)$ which are described in terms of inequalities on parameters. Moreover,
\begin{itemize}
       \item[{\rm(i)}] classes $1$--$18$ have no positive
fixed point and every orbit tends to some fixed point, while classes $19$--$33$ have at least one positive fixed point;
       \item[{\rm(ii)}]
       $\det A<0$ for classes $19$--$25$ and every orbit tends to some fixed point when they have a unique positive fixed point, while $\det A>0$ for classes $26$--$33$.
\end{itemize}
\end{lemma}
\begin{proof}
The conclusion (i) follows from Theorems 4.1--4.2 and Proposition 4.2 in \cite{nwx2023a}, and the conclusion (ii) follows from
Remark 4.4 and Corollary 4.2 in \cite{nwx2023a}.
\end{proof}

\section{The global dynamics} \label{sec:3}
In this section, we focus on the global dynamics of the system \eqref{seasonal-system-equal-db}. We still use $\Psi(t,x)$ to denote the solution map for the system \eqref{seasonal-system-equal-db}  for any initial point $x\in \mathbb{R}^3$. Then
the Poincar\'{e} map $\mathcal{P}$ associated with \eqref{seasonal-system-equal-db} can be written as
$$
\mathcal{P}(x)=\Psi({\omega},x)=\Phi_{\varphi \omega}(lx)
$$
with
$l=e^{-\mu(1-\varphi) \omega}$, where $\Phi_{t}(x)$ is the solution map of the Lotka-Volterra competitive system \eqref{PEVL}. We are interested in the dynamics of the discrete-time dynamical system $\left\{\mathcal{P}^{k}\right\}_{k \geq 0}$ in $\mathbb{R}^{3}_+$.

Hereafter, we always assume that
$$r=b\varphi-\mu(1-\varphi)>0.$$
By Lemma \ref{carrying-simplex}, any Poincar\'{e} map $\mathcal{P}$ of \eqref{seasonal-system-equal-db} has a carrying simplex for all  $\varphi \in (0,1]$ and positive $\omega$, $\mu$, $b$, $r$ and $a_{ij}$. We denote  the set of all maps taking $\mathbb{R}_+^3$ into itself by $\mathcal{T}(\mathbb{R}_+^3)$ and the set of all Poincar\'{e} maps associated with systems  \eqref{seasonal-system-equal-db} which have a carrying simplex by
$$
\mathrm{CLVSE}(3):=\Big\{\begin{array}{l}
    \mathcal{P}\in \mathcal{T}(\mathbb{R}_+^3): \mathcal{P}(x)=\Phi_{\varphi \omega}(lx),~\Phi{~ \text{is~the~solution~map~of}~}\eqref{PEVL}  \\
    \noalign{\smallskip}
      \text{with~} 0<\varphi<1,\omega,\mu,b,r,a_{ij}>0, i,j=1,2,3
\end{array}
\Big\}.
$$
Clearly, $\mathrm{CLVSE}(3)$ is a subset of $\mathrm{CLVS}(3)$. Therefore, the same statements in Lemma \ref{lemma:PFP} are true for $\mathrm{CLVSE}(3)$.

\subsection{Uniqueness and nonuniqueness of positive fixed points}

In this subsection, we focus on the (non)uniqueness of the positive fixed points for $\mathcal{P}$.
Let $\rho(t,\rho_0)$ be the solution of the logistic equation
 \begin{equation}\label{sys-logistic}
\frac{d \rho}{dt}=b\rho(1-\rho)
\end{equation}
through the initial value $\rho_0>0$. Let $\mathcal{M}$ be the Poincar\'e map of the logistic equation with seasonal succession
\begin{equation}\label{logistic-season}
 \begin{dcases}
\frac{d y}{d t}=-\mu y, ~ t\in[k \omega,  k \omega+(1-\varphi) \omega), \\
\frac{d y}{d t}=by\left(1-y\right),~ t\in [k \omega+(1-\varphi) \omega,  (k+1) \omega),
\end{dcases}
\end{equation}
that is,
$$
\mathcal{M}(y)=\rho(\varphi\omega,ly),~~ y\in \mathbb{R}_+.
$$
It is easy to check that $\mathcal{M}$ has a positive fixed point
\begin{equation}
    \rho^*=\frac{1-e^{\mu(1-\varphi) \omega-b \varphi \omega}}{1-e^{-b \varphi \omega}}
\end{equation}
which is globally asymptotically stable if and only if $r>0$ (see \cite{Hsu-Zhao}).

\begin{proposition}\label{P-Map-DF}
For all $x\in \mathbb{R}^3_+$, one has
\begin{equation}\label{poincare-map}
\mathcal{P}^{k}(x)=\rho^*\Phi(k\hat{\rho},\frac{x}{\rho^*}),\quad k=1,2\ldots,
\end{equation}
where
\begin{equation}\label{hat-rho} \hat{\rho}:=\int_0^{\varphi\omega}\rho(s,l\rho^*)ds=\frac{r}{b}\omega.
\end{equation}
\end{proposition}
\begin{proof}
For the solution map $\Phi(t,x)$ of system \eqref{PEVL}, according to \cite{JN-Siam} one has that
$$
\Phi(t,lx)=\rho(t,l\rho^*)\Phi(\int_0^t\rho(s,l\rho^*)ds,\frac{x}{\rho^*}).
$$
By the definition of $\rho^*$,  one has
$$
\begin{aligned}
    \mathcal{P}(x)&=\Phi_{\varphi\omega}(lx)\\    &=\rho(\varphi\omega,l\rho^*)\Phi(\int_0^{\varphi\omega}\rho(s,l\rho^*)ds,\frac{x}{\rho^*})\\   &=\rho^*\Phi(\hat{\rho},\frac{x}{\rho^*}),
\end{aligned}
$$
and moreover,
$$
\mathcal{P}^2(x)=\rho^*\Phi(\hat{\rho},\Phi(\hat{\rho},\frac{x}{\rho^*}))=\rho^*\Phi(2\hat{\rho},\frac{x}{\rho^*}).
$$
It then follows that
$$
\mathcal{P}^{k}(x)=\rho^*\Phi(k\hat{\rho},\frac{x}{\rho^*}),\quad k=1,2\ldots.
$$

It remains to prove \eqref{hat-rho}. Let $z(t)=\rho(t,l\rho^*)$. Then $z(\varphi\omega)=\rho^*$ and $z(t)>0$, $\forall t\geq 0$. Since $z(t)$ satisfies \eqref{sys-logistic}, one has
$$
\frac{\dot{z}(t)}{z(t)}=b(1-z(t)).
$$
Integrating the above equations for $t$ from 0 to $\varphi\omega$, one obtains that
$$
b(\varphi\omega-\hat{\rho})=\mu(1-\varphi)\omega,
$$
which implies
$$
\hat{\rho}=\int_0^{\varphi\omega}\rho(s,l\rho^*)ds=\frac{r}{b}\omega.
$$
Thus, the proof is completed.
\end{proof}

\begin{remark}
Note that Proposition \ref{P-Map-DF} establishes a topological conjugacy between the Poincar\'e map $\mathcal{P}$ and the time-$\hat{\rho}$ map of the flow $\Phi$. Specifically, by introducing the homeomorphism $h: \mathbb{R}^3_+ \to \mathbb{R}^3_+$ defined by $h(x) = \rho^* x$, it follows from Proposition \ref{P-Map-DF} that$$\mathcal{P} = h \circ \Phi(\hat{\rho}, \cdot) \circ h^{-1}.$$This conjugacy ensures that the discrete-time dynamics governed by $\mathcal{P}$ are qualitatively equivalent to the continuous-time evolution of the Lotka--Volterra system \eqref{PEVL} under the time rescaling $\hat{\rho}$. Consequently, the global dynamics for $\mathcal{P}$ can be rigorously analyzed by that of the continuous system \eqref{PEVL}.
\end{remark}

By Hirsch's carrying simplex theory in \cite{hirsch1988}, it is known that the continuous-time dynamical system $\Phi$ has a carrying simplex $\Sigma_{\Phi}$ such that every nontrivial trajectory is asymptotic to the one in $\Sigma_{\Phi}$. By Proposition \ref{P-Map-DF}, we obtain the following result.
\begin{proposition}\label{simplex-relation}
Suppose that $\Lambda$ is an invariant set of $\Phi$. Then
$$
\rho^*\Lambda:=\{\rho^*x:x\in \Lambda\}
$$
is an invariant set of $\mathcal{P}$. In particular,
$$
\Sigma_{\mathcal{P}}=\rho^*\Sigma_{\Phi}
$$
for $\mathcal{P}\in \mathrm{CLVSE}(3)$.
\end{proposition}
\begin{proof}
The conclusions follow from Proposition \ref{P-Map-DF} and the properties of the carrying simplex immediately.
\end{proof}
\begin{proposition}\label{prop:fp-P}
A point $x^*\in \mathbb{R}_+^3$ is a fixed point of $\mathcal{P}$ if and only if
$$
\Phi(\hat{\rho},\frac{x^*}{\rho^*})=\frac{x^*}{\rho^*},
$$
that is,
either $\frac{x^*}{\rho^*}$ is an equilibrium of $\Phi$ or $\Phi(t,\frac{x^*}{\rho^*})$ is a periodic orbit with a period $\hat{\rho}$.
\end{proposition}
\begin{proof}
By Proposition \ref{P-Map-DF},
$$
\mathcal{P}(x)=\rho^*\Phi(\hat{\rho},\frac{x}{\rho^*}).
$$
Therefore,
$$
\mathcal{P}(x^*)=x^* \Leftrightarrow \Phi(\hat{\rho},\frac{x^*}{\rho^*})=\frac{x^*}{\rho^*}.
$$
Then the conclusion is immediate.
\end{proof}
\begin{theorem}
\label{coro:invariant-curve}
Assume that there is a periodic orbit $\Gamma$ with a  minimal period $T_\Gamma>0$ for $\Phi$. Let
$$\eta \triangleq \left(\frac{\omega}{T_\Gamma}:\frac{b}{r}\right).$$
Then $\rho^*\Gamma$ is an invariant closed curve for $\mathcal{P}$, and moreover,
\begin{itemize}
       \item[{\rm(i)}] $\rho^*x$ is a fixed point of $\mathcal{P}$ for all $x\in \Gamma$ if $\eta$ is an integer.
       \item[{\rm(ii)}] Every orbit on $\rho^*\Gamma$ is a periodic orbit of $\mathcal{P}$ with the same minimal period $q$, if $\eta$ is a non-integer rational number such that $\eta = \frac{p}{q}$, where $p, q \in \mathbb{Z}_+$ are coprime and $q > 1$.
        \item[{\rm(iii)}] Every orbit on $\rho^*\Gamma$ is dense with $\rho^*\Gamma$ as the $\omega$-limit set for  $\mathcal{P}$ if $\eta$ is an irrational number.
\end{itemize}
\end{theorem}
\begin{proof}
Since the periodic orbit $\Gamma$ is an invariant set under $\Phi$, it follows from Proposition \ref{simplex-relation} that the curve $\rho^* \Gamma$ is an invariant closed curve for the map $\mathcal{P}$.

We now consider the restriction of $\mathcal{P}$ to the smooth curve $\rho^* \Gamma$, which is topologically conjugate to the circle rotation $R_\eta: \theta \mapsto (\theta + \eta) \pmod 1$ on $S^1 \cong \mathbb{R}/\mathbb{Z}$. The rotation number $\eta$ is defined as the ratio of the map's time-step $\hat{\rho}$ to the minimal period $T_\Gamma$, i.e.,
$$\eta = \frac{\hat{\rho}}{T_\Gamma} = \frac{r \omega}{b T_\Gamma}.$$

(i) If $\eta$ is an integer, then $\hat{\rho} = \eta T_\Gamma$ is a period of $\Gamma$. By Proposition \ref{prop:fp-P}, every point on $\rho^* \Gamma$ is a fixed point of $\mathcal{P}$.

(ii) If $\eta = \frac{p}{q}$ is a non-integer rational number, where $p, q \in \mathbb{Z}_+$ are coprime and $q > 1$, then $q \hat{\rho} = p T_\Gamma$. For any $y = \rho^* x \in \rho^* \Gamma$, the $q$-th iteration of the map satisfies
$$\mathcal{P}^q(\rho^* x) = \rho^* \Phi(q \hat{\rho}, x) = \rho^* \Phi(p T_\Gamma, x) = \rho^* x.$$
Since $q$ is the smallest positive integer such that $q\eta \in \mathbb{Z}$, every orbit on the curve $\rho^* \Gamma$ is a periodic orbit of $\mathcal{P}$ with minimal period $q$.

(iii)  If $\eta$ is an irrational number, then the rotation $R_\eta$ is an irrational rotation. It is a classical result (see, e.g., \cite[Proposition 1.3.3]{Katok1995}) that every orbit of an irrational rotation is dense in $S^1$. Consequently, for any $y \in \rho^* \Gamma$, the orbit $\{\mathcal{P}^n(y)\}_{n\geq 0}$ is dense in $\rho^* \Gamma$, and hence the $\omega$-limit set of every orbit on $\rho^* \Gamma$ is the entire curve $\rho^* \Gamma$.

\end{proof}

\begin{remark}
    Theorem \ref{coro:invariant-curve} implies that the dynamical behavior of system \eqref{seasonal-system-equal-db} depends on the ratio of the frequency of external forced oscillations to the frequency of internal oscillations. The cases
    (i), (ii) and (iii) in Theorem \ref{coro:invariant-curve} correspond to the existence of infinitely many positive
harmonic solutions, subharmonic solutions and quasiperiodic solutions, respectively.
\end{remark}

We are now ready to give the proof of Theorem 1.1.

{\it Proof of Theorem 1.1}.
Recall that system \eqref{PEVL} possesses a carrying simplex. Consequently, any positive periodic orbit of $\Phi_t$ (if it exists) must lie on this simplex, and the Poincar\'e–Bendixson theorem implies the existence of a positive equilibrium. Conclusion (a) therefore follows directly from Proposition \ref{prop:fp-P}.

   Since $\det A \neq 0$, system \eqref{PEVL} can have at most one positive equilibrium. Using Proposition \ref{prop:fp-P}, Theorem \ref{coro:invariant-curve}, and the positive invariance of $\dot{\mathbb{R}}_{+}^{3}$ under $\Phi$, we immediately obtain conclusion (b).

   By Proposition \ref{prop:fp-P} and the positive invariance of $\dot{\mathbb{R}}_{+}^{3}$ under $\Phi_t$, $\mathcal{P}$ has a  positive fixed point $x^*$ if and only if $\frac{x^*}{\rho^*}$ is a positive equilibrium of $\Phi_t$, or  $\Phi(t,\frac{x^*}{\rho^*})$ is a positive  periodic orbit with a period $\hat{\rho}=\frac{r}{b}\omega$. Hence $x^*$ is the unique
 positive fixed point of $\mathcal{P}$ precisely when
$\frac{x^*}{\rho^*}$ is the unique positive equilibrium of $\Phi_t$ and, in addition, either $\Phi_t$ possesses no positive periodic orbit or for any positive periodic orbit $\Gamma$ of $\Phi_t$, $\hat{\rho}$ is not a period of $\Gamma$. This establishes conclusion (c).

\subsection{Global dynamics}
In this subsection, we aim to provide a delicate global dynamical description for the model \eqref{seasonal-system-equal-db}. We define
\begin{equation}\label{beta_alpha}
\alpha_i:=a_{i+1,i+1}-a_{i,i+1},\ \ \beta_i:=a_{i,i-1}-a_{i-1,i-1}, \quad i \ \ {\rm mod}\ \ 3,\ {\rm and}
\end{equation}
\vspace{-4mm}
\begin{equation}\label{Def_zeta}
\begin{array}{rl}
  \zeta:=\beta_1\beta_2\beta_3-\alpha_1\alpha_2\alpha_3. \\
\end{array}
\end{equation}
Set
\begin{equation}
    f_{i}(x)=x_i(b-\sum_{j=1}^{3} a_{i j} x_{j}),\quad i=1,2, 3.
\end{equation}
We first recall the conclusions in \cite{JN-Siam} on the 3-dimensional Lotka-Volterra competition model \eqref{PEVL}.

\begin{lemma}\label{lemma:eigenvalues}
Assume that $\det A \neq 0$ and $x^{*}=(x^{*}_1,x^{*}_2,x^{*}_3)$ is a positive equilibrium of $\Phi$. Then $-b$ is an eigenvalue of $D f(x^{*})$ with associated eigenvector $(x^{*}_1,x^{*}_2,x^{*}_3)$, and the other two eigenvalues, say $\lambda_{1}, \lambda_{2}$, satisfy that
$$
\lambda_{1}+\lambda_{2}=\frac{b}{\operatorname{det} A} \zeta \qquad
{\rm and } \qquad
\lambda_{1} \lambda_{2}\det A>0. $$
\end{lemma}
\begin{proof}
The conclusion follows from Proposition 4.5 and the text below its proof in \cite{JN-Siam}.
\end{proof}

\begin{lemma}\label{lemma:dynamics-PEVL}
Assume $\det A \neq 0$. Then one has the following conclusions.
\begin{itemize}
       \item[{\rm(i)}] There is at most one positive equilibrium for $\Phi$.
        \item[{\rm(ii)}] $\Phi$ has no periodic orbit if there is no positive equilibrium or $\det A<0$ or $\zeta \neq 0$.
        \item[{\rm(iii)}] Suppose that $\Phi$ has only finitely many equilibria. Then every trajectory of $\Phi$ converges to an equilibrium if there is no positive equilibrium, or $\det A<0$, or $\zeta< 0$, or $\zeta> 0$ such that $\partial\Sigma_\Phi$ is not a heteroclinic cycle.
        \item[{\rm(iv)}] Suppose that $\partial\Sigma_\Phi$ is a heteroclinic cycle and there is a positive equilibrium $x^*$. Then every trajectory in $\dot{\Sigma}_\Phi$ except $x^*$  converges to the heteroclinic cycle $\partial\Sigma_\Phi$ if $\zeta>0$.
        \item[{\rm(v)}]  There is a continuum of nontrivial periodic orbits surrounding a center if and only if $\Phi$ has a positive equilibrium with $\zeta=0$ and $\det A>0$; in particular, no periodic orbit is isolated.
\end{itemize}
\end{lemma}

\begin{proof}

Conclusion (i) follows directly from the algebraic form of system \eqref{PEVL}.

Regarding conclusions (ii) and (v), according to Theorem 4.8 and the proof of Theorem 4.4 in \cite{JN-Siam}, the restriction of the flow $\Phi$ to its carrying simplex $\Sigma_\Phi$ is topologically equivalent to a two-dimensional Lotka--Volterra system. It is established that $\Phi$ possesses nontrivial periodic orbits if and only if there exists a positive equilibrium $x^*$ with $\zeta = 0$. In this case, Lemma \ref{lemma:eigenvalues} implies that the two eigenvalues $\lambda_1, \lambda_2$ of $Df(x^*)$, excluding the eigenvalue $-b$, form a pair of purely imaginary eigenvalues, which necessitates $\det A > 0$. Furthermore, according to Theorem 4.2.1 in Hofbauer and Sigmund~\cite{Sigmund}, the existence of a Dulac function for the two-dimensional Lotka–Volterra system rules out isolated periodic orbits. This ensures that any periodic orbit belongs to a continuum surrounding the positive equilibrium, thereby proving (ii) and (v).

For conclusion (iv), based on Theorem 4.12 and the parameter conditions specified in Table 2 of \cite{JN-Siam}, if $\partial\Sigma_\Phi$ is a heteroclinic cycle and a positive equilibrium $x^*$ exists with $\zeta > 0$, then the heteroclinic cycle $\partial\Sigma_\Phi$ attracts every trajectory in $\dot{\Sigma}_\Phi \setminus \{x^*\}$. Finally, conclusion (iii) also follows directly from Theorem 4.12 and Table 2 in \cite{JN-Siam}.
\end{proof}

\begin{proposition}\label{Prop:index}
Assume that $x^{*}$ is an equilibrium of $\Phi$. Then
$$D \mathcal{P}(\rho^{*} x^{*})=\exp \left\{D f(x^{*}) \hat{\rho}\right\},$$
and moreover, if $\det A \neq 0$ then
$$
\operatorname{ind}(\mathcal{P},\rho^{*} x^{*})=\operatorname{sgn} \det A,
$$
where $\operatorname{ind}(\mathcal{P},\rho^{*} x^{*})$ is
the fixed point index of $\mathcal{P}$ at its fixed point $\rho^{*} x^{*}$.
\end{proposition}
\begin{proof}
By Proposition \ref{P-Map-DF}, $\mathcal{P}(x)=\rho^{*} \Phi(\hat{\rho}, \frac{x}{\rho^{*}})$,
so
$D \mathcal{P}(\rho^{*} x^{*})=D_{x} \Phi(\hat{\rho}, x^{*})$. Let $W(t,x)=D_{x} \Phi(t,x)$ and $U(t,x)=D f(\Phi(t,x))$. Then
\begin{equation}\label{equ-DPhi}
    \frac{d W(t,x)}{d t}=U(t,x) \cdot W(t,x), \quad W(0,x)=I.
\end{equation}
Since $\Phi(t, x^{*})=x^*$ for all $t\in \mathbb{R}$,
one has
$$D_{x} \Phi(t, x^{*})=\exp \left\{D f(x^{*})t\right\},$$
which implies that
$$D \mathcal{P}(\rho^{*} x^{*})=\exp \left\{D f(x^{*}) \hat{\rho}\right\}.$$
Note that each eigenvalue of
$D \mathcal{P}(\rho^{*} x^{*})$ is the exponential of an eigenvalue of $D f(x^{*}) \hat{\rho}$, so Lemma \ref{lemma:eigenvalues} implies that $D \mathcal{P}(\rho^{*} x^{*})$ has only one eigenvalue greater than one if $\det A<0$ while it has zero or two eigenvalues greater than one if $\det A>0$. Therefore, if $\det A\neq 0$, then
$
\operatorname{ind}(\mathcal{P},\rho^{*} x^{*})=\operatorname{sgn} \det A
$.  \end{proof}

\begin{proposition}\label{trivial-dynamics}
Let $\mathcal{P}\in \mathrm{CLVSE}(3)$. Assume that $\mathrm{Fix}(\mathcal{P}, \dot{\mathbb{R}}^3_+)=\emptyset$ or $\mathrm{Fix}(\mathcal{P},$ $\dot{\mathbb{R}}^3_+)=\{q\}$ with $\mathrm{ind}(\mathcal{P},q)=-1$. Then $\mathcal{P}$ has trivial dynamics, i.e., the omega limit set of any orbit only contains fixed points. Furthermore, if $\mathcal{P}$ has only finitely many fixed points, then every orbit converges to some fixed point.
\end{proposition}
\begin{proof}
The conclusion can be proved in a manner similar to Proposition 4.1 in \cite{nwx2023a}.
\end{proof}


\begin{proposition}\label{thm:dy-P}
For $\mathcal{P}\in \mathrm{CLVSE}(3)$, the following conclusions hold:
\begin{itemize}
       \item[{\rm(i)}] $\mathcal{P}$ has trivial dynamics if the linear algebraic system $(A x^{\tau})_{i}=b$, $i=1,2,3$ has no positive solution.
       \item[{\rm(ii)}] $\mathcal{P}$ has trivial dynamics if $\det A<0$, and moreover, $\mathcal{P}$ has at most one positive  fixed point which is a saddle on $\Sigma_{\mathcal{P}}$ if it exists.
       \item[{\rm(iii)}] Suppose that $\det A>0$. If $\zeta< 0$, then $\mathcal{P}$ has trivial dynamics,  and moreover, $\mathcal{P}$ has at most one positive  fixed point which is an attractor on $\Sigma_{\mathcal{P}}$ if it exists.
       \item[{\rm(iv)}] Suppose that $\det A>0$. If $\zeta> 0$ and $\partial \Sigma_\mathcal{P}$ is not a  heteroclinic cycle, then $\mathcal{P}$ has trivial dynamics,  and moreover, $\mathcal{P}$ has at most one positive  fixed point which is a repeller on $\Sigma_{\mathcal{P}}$ if it exists.
       \item[{\rm(v)}] Suppose that $\det A>0$ and $\mathcal{P}$ has a positive fixed point $p$. If $\zeta> 0$ and $\partial\Sigma_\mathcal{P}$ is a  heteroclinic cycle, then every trajectory in $\dot{\Sigma}_\mathcal{P}$ except $p$ converges to the heteroclinic cycle $\partial \Sigma_\mathcal{P}$.
       \item[{\rm(vi)}] If $\det A > 0$, $\zeta = 0$, and $(A x^{\tau})_i = b$ ($i=1,2,3$) possesses a unique positive solution, then $\mathcal{P}$ admits a continuum of invariant closed curves surrounding a positive fixed point; in particular, no such curve is isolated. Furthermore, each curve is characterized by one of the following: every point is a fixed point, every orbit is periodic, or every orbit is dense.
\end{itemize}
\end{proposition}
\begin{proof}
(i) If $(A x^{\tau})_{i}=b$, $i=1,2,3$ has no positive solution, then $\Phi$ does not have positive equilibrium and periodic orbit by Lemma \ref{lemma:dynamics-PEVL}. Therefore, $\mathcal{P}$ has no positive fixed point by Proposition \ref{prop:fp-P}. Then the conclusion (i) follows from Proposition \ref{trivial-dynamics}.

(ii) If $\det A<0$, then $\Phi$ has at most one positive equilibrium and no periodic orbit by Lemma  \ref{lemma:dynamics-PEVL}, and hence $\mathcal{P}$ has at most one positive fixed point by Proposition \ref{prop:fp-P}.
If $\mathcal{P}$ has a positive fixed point, say $p$, then
$$
\operatorname{ind}(\mathcal{P},p)=\operatorname{sgn} \det A=-1
$$
by Proposition \ref{Prop:index}.
Thus, one has either
$\mathrm{Fix}(\mathcal{P}, \dot{\mathbb{R}}^3_+)=\emptyset$
 or $\mathrm{Fix}(\mathcal{P}, \dot{\mathbb{R}}^3_+)=\{p\}$ with
$\mathrm{ind}(\mathcal{P},p)=-1$, and hence the conclusion is true by Proposition \ref{trivial-dynamics}.

(iii) If $\det A>0$ and $\zeta< 0$, then $\mathcal{P}$ has at most one positive fixed point by Proposition \ref{prop:fp-P} and Lemma  \ref{lemma:dynamics-PEVL}(i)-(ii). Therefore, $\mathcal{P}$ has only finitely many fixed points, which implies that $\Phi$ has finitely many fixed points by Proposition \ref{prop:fp-P}, and moreover, it follows from Lemma  \ref{lemma:dynamics-PEVL}(iii) and \eqref{poincare-map} that every trajectory of $\mathcal{P}$ converges to a fixed point. Assume that $\mathcal{P}$ has a positive fixed point, say $p$. By Lemma \ref{lemma:eigenvalues} and Proposition  \ref{Prop:index}, $e^{-b\hat{\rho}}$ is an eigenvalue of $D\mathcal{P}(p)$ and $D\mathcal{P}(p)$ has a positive eigenvector associated with the eigenvalue $e^{-b\hat{\rho}}$. Let $\lambda_1,\lambda_2$ be the other two eigenvalues of $D\mathcal{P}(p)$ besides $e^{-b\hat{\rho}}$. Then Lemma \ref{lemma:eigenvalues} and Proposition \ref{Prop:index} imply that each of $\lambda_1$ and $\lambda_2$ has modulus less than one if $\zeta<0$ and $\det A>0$ while has modulus greater than one if $\zeta>0$ and $\det A>0$. Thus, $p$ is a hyperbolic fixed point of $\mathcal{P}$. Since $(D\mathcal{P}(p))^{-1}$ is a positive matrix (see Lemma \ref{prop-lambda}), $e^{-b\hat{\rho}}$ is the eigenvalue of smallest modulus by Perron-Frobenius theorem.
By \cite[Theorem 4.6]{MNR2019}, the local dynamics of  $p$ is determined by $\lambda_1,\lambda_2$, and hence $p$ is an attractor (a repeller) if $\zeta<0$ ($\zeta>0$).

(iv) The result can be obtained by the similar arguments as (iii).

(v) Proposition \ref{simplex-relation} and Lemma  \ref{lemma:dynamics-PEVL}(ii) implies that $\partial\Sigma_\Phi$ is a heteroclinic cycle and $\Phi$ has a positive equilibrium if $\partial\Sigma_\mathcal{P}$ is a  heteroclinic cycle and $\mathcal{P}$ has a positive fixed point. Then the result follows from Lemma \ref{lemma:dynamics-PEVL}(iv) and \eqref{poincare-map}.

(vi) If $\det A>0$, $\zeta=0$, and $(A x^{\tau})_{i}=b, i=1,2,3$ has a unique positive solution, then $\Phi$ has a positive equilibrium with $\zeta=0$ and $\det A>0$. By Lemma \ref{lemma:dynamics-PEVL}(v), $\Phi$ possesses a continuum of nontrivial periodic orbits surrounding a center, where each periodic orbit is non-isolated. By invoking Proposition \ref{prop:fp-P} and Theorem \ref{coro:invariant-curve}, it follows that $\mathcal{P}$ possesses a continuum of invariant closed curves surrounding a positive fixed point, and each invariant closed curve is non-isolated. The assertion regarding the dynamics on these curves, including whether they consist of fixed points, periodic orbits, or dense orbits, follows directly from the properties established in Theorem \ref{coro:invariant-curve}.

\end{proof}

Since the dynamics of classes 1--18 for  $\mathrm{CLVSE}(3)$ is clear by Lemma \ref{lemma:PFP}, we only focus on the classes 19--33 here. Note that classes 19--33 have finitely many fixed points on $\partial \mathbb{R}_+^3$.

\begin{theorem}[Global dynamics for classes 19--25]\label{thm:whole-dynamics}
Any map from classes 19--25 for $\mathrm{CLVSE}(3)$ has a unique positive fixed point, and moreover, every orbit converges to some fixed point in these classes.  The phase portraits on the carrying simplices for these classes are given in Table \ref{biao0}.
\end{theorem}
\begin{proof}
Since $\det A<0$ for classes 19--25 in $\mathrm{CLVSE}(3)$ by Lemma \ref{lemma:PFP}, it follows from Proposition \ref{thm:dy-P}(ii) and Lemma \ref{lemma:PFP} that there is a unique positive fixed point, say $p$, for any map in these classes, and moreover, every orbit converges to some fixed point. By Lemma \ref{lemma:eigenvalues} and Proposition \ref{Prop:index},  $1$ is not an eigenvalue of $D\mathcal{P}(p)$. Then together with Lemma \ref{prop-lambda} we know that all three eigenvalues of $D\mathcal{P}(p)$, say $\lambda, \lambda_1,\lambda_2$, are positive real numbers with
$$0<\lambda<\lambda_1<1<\lambda_2.$$
Therefore, $p$ is a hyperbolic saddle whose stable manifold and unstable manifold on the carrying simplex are simple curves, and the whole dynamics on the carrying simplices is as shown in Table \ref{biao0} (see \cite[Corollary 5.4]{MNR2019}).
\end{proof}

\begin{figure}[h]
	\begin{center}
		\includegraphics[width=0.25\textwidth]{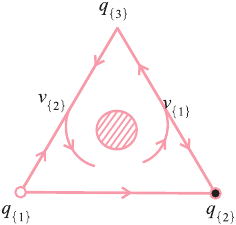}
		\caption{The phase portrait on the carrying simplex for class $26$. A fixed point is represented by a closed dot $\bullet$ if it attracts on the carrying simplex, by an open dot $\circ$ if it repels, and by the intersection of its stable and unstable manifolds if it is a saddle. The circle hatched with oblique lines
 denotes a region of unknown dynamics where there might be more than one fixed points or other complex dynamics such as invariant closed curves.} \label{fig-26}
	\end{center}
\end{figure}

\begin{lemma}[\cite{nwx2023a}]\label{clvs-26}
The map $\mathcal{P}\in \mathrm{CLVS}(3)$ is in class $26$ if and only if there is a permutation $\sigma$ of the indices $\{1, 2 , 3\}$, after which it {\rm(}still denoted by $\mathcal{P}${\rm)} satisfies the following inequalities
\begin{itemize}
  \item[{\rm(i)}] $\gamma_{12}>0, \gamma_{13}>0, \gamma_{21}<0, \gamma_{23}<0, \gamma_{31}>0, \gamma_{32}<0$;
  \item[{\rm(ii)}] $a_{12}\beta_{23}+a_{13}\beta_{32}>r_1$;
  \item[{\rm(iii)}] $a_{21}\beta_{13}+a_{23}\beta_{31}<r_2$.
\end{itemize}
where
\begin{equation}\label{gamma-beta}
    \gamma_{ij}=a_{ii}r_j-a_{ji}r_i,\quad \beta_{ij}=\frac{a_{jj}r_i-a_{ij}r_j}{a_{ii}a_{jj}-a_{ij}a_{ji}},\quad i,j=1,2,3,  ~i\neq j.
\end{equation}
In this case, there are two planar fixed points $v_{\{1\}},v_{\{2\}}$ which are saddles and at least one positive fixed point. The phase portrait on the carrying simplex is as shown in Fig. \ref{fig-26}.
\end{lemma}
\vspace{-0.2cm}

\begin{theorem}[Global dynamics for class 26]\label{prop:class_26}
Within the class $26$ for $\mathrm{CLVSE}(3)$, there are totally three dynamical scenarios depending upon the sign of $\zeta$.
\begin{enumerate}[{\rm (i)}]
\item Assume $\zeta\neq 0$ for the map $\mathcal{P}$ in class $26$. Then $\mathcal{P}$ has a unique positive fixed point $p$, and $p$ is an attractor if $\zeta< 0$ while it is a repeller if $\zeta>0$. Moreover, every orbit converges to some fixed point.

\item Assume $\zeta=0$ for the map $\mathcal{P}$ in class $26$. Then there is a continuum of invariant closed curves enclosed by a heteroclinic cycle, and moreover, on each curve either every point is a fixed point, or every orbit is periodic, or every orbit is dense on the curve.
\end{enumerate}
The phase portraits are presented in {\rm Fig. \ref{fig:class26}}.
\begin{figure}[h]
    \centering
    \begin{tabular}{ccc}
        \subfigure[$\zeta<0$]{
            \label{26_a}
            \begin{minipage}[b]{0.25\textwidth}
                \centering
                \includegraphics[width=\textwidth]{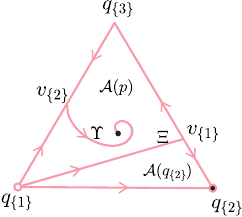}
            \end{minipage}
        } &
        \subfigure[$\zeta>0$]{
            \label{26_b}
            \begin{minipage}[b]{0.25\textwidth}
                \centering
                \includegraphics[width=\textwidth]{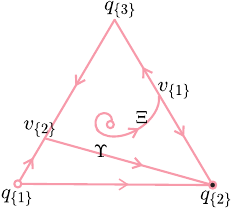}
            \end{minipage}
        }
        \subfigure[$\zeta=0$]{
            \label{26_c}
            \begin{minipage}[b]{0.25\textwidth}
                \centering                \includegraphics[width=\textwidth]{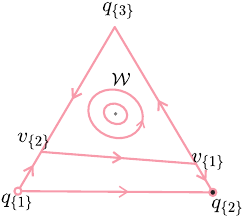}
            \end{minipage}
        } & \\
    \end{tabular}
    \caption{The phase portraits for the three cases in class
    $26$. (a) $\mathcal{A}(p)$ (resp. $\mathcal{A}(q_{\{2\}})$) is the  basin of attraction in the interior of $\Sigma_\mathcal{P}$ for the attractor $p$ (resp. $q_{\{2\}}$). (b) All orbits in the interior of $\Sigma_\mathcal{P}$ except the positive fixed point and the stable manifold of $v_{\{1\}}$ approach the attractor $q_{\{2\}}$. (c) The region $\mathcal{W}$ inside $\Sigma_\mathcal{P}$ is filled with a continuum of invariant closed curves surrounding a positive fixed point; the boundary of  $\mathcal{W}$ is a heteroclinic cycle. The $\Xi$ (resp. $\Upsilon$) is the stable (resp. unstable) manifold of $v_{\{1\}}$ (resp. $v_{\{2\}}$) in $\Sigma_{\mathcal{P}}$. The fixed point notation is as in Fig. \ref{fig-26}.}
    \label{fig:class26}
\end{figure}
\end{theorem}

\begin{proof}
Since $\det A>0$ for the class 26 in $\mathrm{CLVSE}(3)$ by Lemma \ref{lemma:PFP} and the boundary of the carrying simplex is not a heteroclinic cycle, there is a unique positive fixed point for any map in class 26 by Proposition \ref{thm:dy-P} and Lemma \ref{lemma:PFP} if $\zeta\neq 0$, and moreover, the stability of $p$ follows from Proposition \ref{thm:dy-P}. The stability of all boundary fixed points follows from Lemma \ref{clvs-26}.  Denote by $\Xi$ the stable manifold of $v_{\{1\}}$ in $\Sigma_{\mathcal{P}}$, and by $\Upsilon$ the unstable manifold of $v_{\{2\}}$ in $\Sigma_{\mathcal{P}}$, respectively.

For the dynamics of the maps in class 26, it suffices to study that of the 3-dimensional Lotka-Volterra competition model \eqref{PEVL} by Proposition \ref{P-Map-DF}. Specifically, $\Phi$ has a carrying simplex
$$
\Sigma_{\Phi}=\frac{1}{\rho^*}\Sigma_{\mathcal{P}},
$$
 and three axial fixed points $\hat{q}_{\{i\}}=\frac{1}{\rho^*}q_{\{i\}}$, ($i=1,2,3$), two planar fixed points $\hat{v}_{\{i\}}=\frac{1}{\rho^*}{v}_{\{i\}}$, ($i=1,2$), and a unique positive fixed point $\hat{p}=\frac{1}{\rho^*}p$, where $\mathcal{P}$ is in the class 26. Moreover, Proposition \ref{Prop:index} implies that $\hat{q}_{\{i\}}$ and $q_{\{i\}}$, $\hat{v}_{\{i\}}$ and ${v}_{\{i\}}$ have the same local dynamics, respectively. Therefore,
$\hat{v}_{\{1\}}$ and $\hat{v}_{\{2\}}$ are saddles, and the unstable manifold $\hat{\Upsilon}$  of $\hat{v}_{\{2\}}$ in $\Sigma_{\Phi}$ can not coincide with the stable manifold $\hat{\Xi}$ of $\hat{v}_{\{1\}}$ in $\Sigma_{\Phi}$ if $\zeta\neq 0$ because by \cite{JN-Siam} there is no heteroclinic cycle for $\Phi$  if $\zeta\neq 0$ when the parameters are in the class 26, where $\hat{\Upsilon}=\frac{1}{\rho^*}\Upsilon$ and $\hat{\Xi}=\frac{1}{\rho^*}\Xi$. By Lemma \ref{lemma:dynamics-PEVL}(iii),  every trajectory of $\Phi$ converges to an equilibrium if $\zeta\neq 0$. It follows from the Poincar\'{e}-Bendixson Theorem that the unstable manifold $\hat{\Upsilon}$ tends towards the attractor $\hat{p}$ and the stable manifold $\hat{\Xi}$ flows away from the repeller $\hat{q}_{\{1\}}$ if $\zeta<0$, and moreover, the manifold $\hat{\Xi}$ divides $\dot{\Sigma}_\Phi\setminus \hat{\Xi}$ into two parts, corresponding to the basins of attraction in $\dot{\Sigma}_\Phi$ of the two attractors, $\hat{p}$ and $\hat{q}_{\{2\}}$, respectively. If $\zeta>0$, the unstable manifold $\hat{\Upsilon}$ is towards the attractor $\hat{q}_{\{2\}}$ and the stable manifold $\hat{\Xi}$ of $\hat{v}_{\{1\}}$ flows away from the repeller $\hat{p}$. Trajectories in $\dot{\Sigma}_\Phi\setminus (\hat{\Xi}\cup \{\hat{p}\})$ approach the attractor $\hat{q}_{\{2\}}$. Then by the above arguments and Proposition \ref{P-Map-DF} we obtain the results of (i) and the phase portraits on the carrying simplex (see Figs. \ref{26_a} and \ref{26_b}).

 For $\zeta=0$, the unstable manifold $\hat{\Upsilon}$ coincides with the stable manifold $\hat{\Xi}$, and hence it is a saddle-connection of $\hat{v}_{\{1\}}$ and $\hat{v}_{\{2\}}$. Therefore, $\hat{q}_{\{3\}}\to \hat{v}_{\{2\}}\to \hat{v}_{\{1\}}\to \hat{q}_{\{3\}}$ forms a heteroclinic cycle, and the area $\hat{\mathcal{W}}$ in $\dot{\Sigma}_\Phi$ enclosed by the heteroclinic cycle is filled with a continuum of nontrivial periodic orbits around a positive equilibrium $\hat{p}$ by Lemma \ref{lemma:dynamics-PEVL}. Trajectories in $\dot{\Sigma}_\Phi\setminus (\hat{\mathcal{W}}\cup \hat{\Xi})$  flow away from $\hat{q}_{\{1\}}$ and approach the attractor $\hat{q}_{\{2\}}$ by Poincar\'{e}-Bendixson Theorem. Then the results of (ii) and the phase portrait on the carrying simplex follow from Proposition \ref{P-Map-DF} and Theorem \ref{coro:invariant-curve} immediately (see Fig. \ref{26_c}).
\end{proof}

\begin{figure}[h]
	\begin{center}		\includegraphics[width=0.25\textwidth]{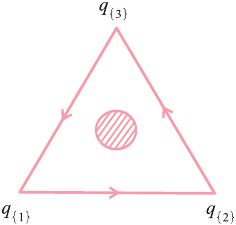}
		\caption{The phase portrait on the carrying simplex for class $27$. The fixed point notation is as in Fig. \ref{fig-26}.} \label{fig-27}
	\end{center}
\end{figure}

\begin{lemma}[\cite{nwx2023a}]
The map $\mathcal{P}\in \mathrm{CLVS}(3)$ is in the class $27$ if and only if there is a permutation $\sigma$ of the indices $\{1, 2 , 3\}$, after which it {\rm (}still denoted by $\mathcal{P}${\rm)} satisfies the following inequalities
$$\gamma_{12}>0, \gamma_{13}<0, \gamma_{21}<0, \gamma_{23}>0, \gamma_{31}>0, \gamma_{32}<0,$$
where
$\gamma_{ij}$ are given by \eqref{gamma-beta}.
In this case, there is at least one positive fixed point and $\partial \Sigma_\mathcal{P}$ is a heteroclinic cycle. The phase portrait on $\Sigma_\mathcal{P}$ is as shown in Fig. \ref{fig-27}.
\end{lemma}

\begin{lemma}\label{lemma:hetero-per}
Assume that the map $\mathcal{P}\in \mathrm{CLVSE}(3)$ is in the class $27$. Then $\partial \Sigma_\mathcal{P}$ is a heteroclinic cycle, and moreover, $\partial \Sigma_\mathcal{P}$ is repelling {\rm(}resp. attracting{\rm )} if $\vartheta>0~({resp.} <0)$, where
\begin{equation}\label{vartheta}
    \vartheta=w_{12} w_{23} w_{31}+w_{21}w_{13}w_{32},
\end{equation}
with $w_{ij}=r-a_{ji}\frac{r}{a_{ii}}$,  $i,j=1,2,3$, $i\neq j$.
\end{lemma}
\begin{proof}
    Since the class $27$ in $\mathrm{CLVSE}(3)$ is a subset of class $27$ in $\mathrm{CLVS}(3)$, one has that $\partial \Sigma_\mathcal{P}$ is a heteroclinic cycle. The remaining results follow from \cite[Corollary 4.1]{nwx2023a}.
\end{proof}

\begin{theorem}[Global dynamics for class 27]\label{prop:class_27}
Within the class $27$ for $\mathrm{CLVSE}(3)$, there are totally three dynamical scenarios depending upon the sign of $\zeta$.
\begin{enumerate}[{\rm (i)}]
\item Assume $\zeta \neq 0$ for the map $\mathcal{P}$ in class $27$. Then $\mathcal{P}$ has a unique positive fixed point $p$. Moreover,  $p$ is globally attracting in $\dot{\Sigma}_\mathcal{P}$ if $\zeta < 0$, while the heteroclinic cycle $\partial \Sigma_\mathcal{P}$ is globally attracting in $\dot{\Sigma}_\mathcal{P}$ if $\zeta>0$.
\item Assume $\zeta=0$  for the map $\mathcal{P}$ in class $27$.  Then $\dot{\Sigma}_\mathcal{P}$ is all covered by a continuum of invariant closed curves around a positive fixed point, and on each curve either every point is a fixed point, or every orbit is periodic, or every orbit is dense on the curve.
\end{enumerate}
The phase portraits are presented in {\rm Fig. \ref{fig:class27}}.
\vspace{-0.15cm}
\begin{figure}[h]
    \centering
    \begin{tabular}{ccc}
        \subfigure[$\zeta<0$]{
            \label{27_a}
            \begin{minipage}[b]{0.23\textwidth}
                \centering
                \includegraphics[width=\textwidth]{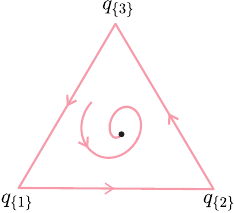}
            \end{minipage}
        } &
        \subfigure[$\zeta>0$]{
            \label{27_b}
            \begin{minipage}[b]{0.23\textwidth}
                \centering
        \includegraphics[width=\textwidth]{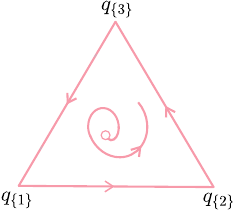}
            \end{minipage}
        } &
        \subfigure[$\zeta=0$]{
            \label{27_c}
            \begin{minipage}[b]{0.23\textwidth}
                \centering
         \includegraphics[width=\textwidth]{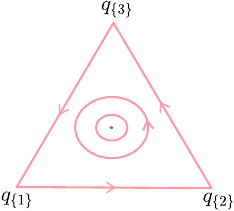}
            \end{minipage}
        }  \\
    \end{tabular}
    \caption{The phase portraits for the three cases in class $27$. (a) The positive fixed point is globally asymptotically stable in $\dot{\mathbb{R}}_+^3$. (b) All orbits in the interior of $\Sigma_\mathcal{P}$ except the positive fixed point approach the heteroclinic cycle. (c) The interior of $\Sigma_\mathcal{P}$ is entirely filled by a continuum of invariant closed curves surrounding a positive fixed point.
 The fixed point notation is as in Fig. \ref{fig-26}.}
    \label{fig:class27}
\end{figure}
\end{theorem}

\begin{proof}
First, $\partial \Sigma_\mathcal{P}$ is a heteroclinic cycle for any map $\mathcal{P}$ in class 27 for $\mathrm{CLVSE}(3)$ by Lemma \ref{lemma:hetero-per}. It is easy to check that $\mathrm{sgn}(\vartheta)=-\mathrm{sgn}( \zeta)$ by \eqref{vartheta}, so $\partial \Sigma_\mathcal{P}$ is a repelling (resp. an attracting) heteroclinic cycle if $\zeta< 0$ (resp. $\zeta> 0$). Since $\det A>0$ for the class 27 in $\mathrm{CLVSE}(3)$ by Lemma \ref{lemma:PFP}, there is a unique positive fixed point for any map in class 27 by Proposition \ref{thm:dy-P} and Lemma \ref{lemma:PFP} if $\zeta\neq 0$. Now the conclusions follow from Proposition \ref{thm:dy-P}(iii), (v) and (vi),  respectively.
 \end{proof}
 \begin{remark}[Multiplicity of positive fixed points]\label{rmk:MFP}
Note that there exists an invariant closed curve $\rho^*\Gamma$ for  $\mathcal{P}$ in the classes 26 and 27 with $\zeta=0$, where $\Gamma$ is a periodic orbit of $\Phi$ with a minimal period $T_\Gamma>0$. By Theorem \ref{coro:invariant-curve}, every point on $\rho^*\Gamma$ is a fixed point if
$\frac{\omega}{T_\Gamma}=k\frac{b}{r}$ with $k>0$ being an integer. In particular, if we set
$$
\omega=\frac{b}{r}T_\Gamma,
$$
then every point on $\rho^*\Gamma$ is a positive fixed point.
Recall that $\det A\neq 0$ for any map in the classes 26 and 27, so $\mathcal{P}$ can have more than one positive fixed point when $\det A\neq 0$, which is a significant difference from the 2-dimensional case.
\end{remark}
\begin{theorem}[Global dynamics for classes 28--33]\label{thm:class_28_33}
Within classes $28$ to $33$, the positive fixed point is unique. Moreover, it is an attractor {\rm(}resp. a repeller{\rm)} on the carrying simplex in classes $29$, $31$, $33$ {\rm(}resp. classes $28$, $30$, $32${\rm)}, and every trajectory converges to some fixed point. The phase portraits on the carrying simplex for classes $28$ to $33$ are given in Table \ref{biao0}.
\end{theorem}
\begin{proof}
It suffices to show that $\zeta\neq 0$ for each map in classes $28$ to $33$. Once this is verified, the stability of $p$ follows directly from Proposition \ref{thm:dy-P}.

(1) For class $28$, parameter condition (i) in  Table \ref{biao0}(28) gives
$$(a_{13}-a_{33})(a_{21}-a_{11})(a_{32}-a_{22})>0$$
and
$$(a_{22}-a_{12})(a_{33}-a_{23})(a_{11}-a_{31})<0,$$
hence
$$\zeta= (a_{13}-a_{33})(a_{21}-a_{11})(a_{32}-a_{22})-(a_{22}-a_{12})(a_{33}-a_{23})(a_{11}-a_{31})>0.$$
For class $32$, the same argument yields $\zeta>0$. Thus, in both cases $p$ is a repeller on the carrying simplex.

(2) For classes $29$ and $33$, the corresponding condition (i) in Table \ref{biao0} implies $\zeta<0$. Therefore, in these two classes $p$ is an attractor on the carrying simplex.

(3) For class $30$, we first show that $a_{12} > a_{32}$. Suppose otherwise; then by
condition (i) in  Table \ref{biao0}(30),
\begin{align*}
     &a_{12}(a_{33}-a_{23})+a_{13}(a_{22}-a_{32})+a_{23}a_{32}-a_{22}a_{33}\\
     \geq& a_{32}(a_{33}-a_{23})+a_{13}(a_{22}-a_{32})+a_{23}a_{32}-a_{22}a_{33}\\
    =&(a_{32}-a_{22})(a_{33}-a_{13})>0,
\end{align*}
contradicting (ii). Using (ii) and (iii) of Table \ref{biao0}(30) we then obtain
$$a_{23}-a_{33}>(a_{32}-a_{22})(a_{33}-a_{13})/(a_{12}-a_{32})$$
 and
$$a_{31}-a_{11}>(a_{21}-a_{11})(a_{12}-a_{32})/(a_{12}-a_{22}).$$
Therefore,
\begin{align*}
     \zeta =& (a_{13}-a_{33})(a_{21}-a_{11})(a_{32}-a_{22})+(a_{12}-a_{22})(a_{23}-a_{33})(a_{31}-a_{11})\\
    >& (a_{13}-a_{33})(a_{21}-a_{11})(a_{32}-a_{22})+(a_{33}-a_{13})(a_{21}-a_{11})(a_{32}-a_{22})\\=&0,
\end{align*}
so in this case $p$ is a repeller on the carrying simplex.

(4) For class $31$, conditions (i)-(iii) in Table \ref{biao0}(31) directly give
$$a_{13}-a_{33}<(a_{23}-a_{33})(a_{12}-a_{32})/(a_{22}-a_{32})$$
 and
 $$a_{31}-a_{11}<(a_{21}-a_{11})(a_{12}-a_{32})/(a_{12}-a_{22}).$$
Consequently,
\begin{align*}
    \zeta =&   (a_{13}-a_{33})(a_{21}-a_{11})(a_{32}-a_{22})+(a_{12}-a_{22})(a_{23}-a_{33})(a_{31}-a_{11})\\
    <& -(a_{23}-a_{33})(a_{21}-a_{11})(a_{12}-a_{32})+(a_{23}-a_{33})(a_{21}-a_{11})(a_{12}-a_{32})\\
    =&0,
    \end{align*}
and hence $p$ is an attractor on the carrying simplex.
\end{proof}




\section{Discussion} \label{sec:7}
In this paper, we establish necessary and sufficient conditions for the uniqueness or nonuniqueness of positive fixed points of the Poincar\'e map associated with the Lotka-Volterra competitive system \eqref{seasonal-system-equal-db} subject to seasonal succession. This system is characterized by an identical intrinsic growth rate during the favorable season and an identical death rate during the unfavorable season. We further prove that the number of positive fixed points is determined by the existence and periods of periodic orbits in the classical Lotka-Volterra competitive system \eqref{PEVL}, which shares the same intrinsic growth rate as system \eqref{seasonal-system-equal-db}.

Furthermore, we delineate all possible global dynamical scenarios of the Poincar\'e map associated with system \eqref{seasonal-system-equal-db}, which amount to a total of 37 distinct dynamical scenarios. Specifically, among the 33 equivalence classes established in \cite{nwx2023a}, dynamical classes 1–18, which possess no positive fixed point, exhibit trivial dynamics. In the present work, we prove that a unique positive fixed point exists in dynamical classes 19–25 and 28–33, and these classes also display trivial dynamics. Each of these classes corresponds to exactly one dynamical scenario.

In contrast, classes 26 and 27 each encompass three distinct dynamical scenarios. Subcases 26(a), 26(b), 27(a), and 27(b) possess a unique positive fixed point, and every orbit converges to a fixed point in subcases 26(a), 26(b), and 27(a). The dynamics in subcases 26(c), 27(b), and 27(c) are much richer and more complex:
\begin{itemize}
    \item Class 26(c) displays the coexistence of an asymptotically stable fixed point, a continuum of invariant closed curves, and a heteroclinic cycle on the carrying simplex;
    \item Class 27(b) features a heteroclinic cycle that attracts all points except a positive fixed point;
    \item Class 27(c) exhibits a continuum of invariant closed curves covering the entire interior of the carrying simplex.
\end{itemize}

Returning to system \eqref{seasonal-system-equal-db}, our results imply that every solution belonging to classes 1–25, 26(a), 26(b), 27(a), and 28–33 converges asymptotically to an $\omega$-periodic solution. In particular, a globally attracting positive $\omega$-periodic solution exists in classes 29, 31, 33, and 27(a).
More notably, our findings demonstrate that the Poincar\'e map of system \eqref{seasonal-system-equal-db} can indeed admit multiple positive fixed points, which corresponds to the existence of multiple positive $\omega$-periodic solutions in the system.

However, it is still unknown when the positive fixed point is (non)unique for the Poincar\'e map of system \eqref{seasonal-system} which has distinct intrinsic growth rates or death rates. An interesting question is whether the number of positive fixed points is related to the periodic orbits of the  Lotka-Volterra competitive system with distinct intrinsic growth rates. Note that for the 3-dimensional system \eqref{seasonal-system-equal-db} with identical intrinsic growth rate and death rate, there are only nonisolated invariant closed curves. Though isolated invariant closed curve has been found in \cite{nwx2023a} for the Poincar\'e map of the 3-dimensional system \eqref{seasonal-system} with distinct intrinsic growth rates and death rates, whether it has nonisolated invariant closed curves is also an interesting problem.

\section*{Acknowledgments}
The authors are greatly indebted to the anonymous referees
for very careful reading and many suggestions which led to many improvements to
the earlier version of this paper. The authors are very grateful to Prof. Sze-Bi Hsu for his valuable and useful discussions and suggestions.

\bibliographystyle{siamplain}
\bibliography{refs}

\appendix
\section{The global dynamics for $\mathrm{CLVSE}(3)$}\label{appendix}
\setcounter{table}{0}
\begin{center}
  \begin{longtable}{c@{\extracolsep{\fill}}c@{\extracolsep{\fill}}c@{\extracolsep{\fill}}}
\caption{The $37$ global dynamical
scenarios in the 33 equivalence classes for $\mathrm{CLVSE}(3)$ which are described in terms of inequalities on parameters of the model, where
$\gamma_{ij}=a_{ii}r-a_{ji}r$, $\beta_{ij}=\frac{a_{jj}r-a_{ij}r}{a_{ii}a_{jj}-a_{ij}a_{ji}}$, $i,j=1,2,3$, $i\neq j$, and $\zeta$ is given by \eqref{Def_zeta}. Each carrying simplex is given by a representative map of that class. The fixed point notation is as in Fig. \ref{fig-26}. } \\[-2pt]
        \hline
          \footnotesize {Classes} &   \footnotesize{Parameter conditions} &   \footnotesize{Phase Portraits} \\
        \hline
        \endfirsthead
        \caption[]{(continued)}\\
        \hline
          \footnotesize{Classes} &   \footnotesize{Parameter conditions} &   \footnotesize{Phase Portraits}  \\
        \hline
&&\\
        \endhead
        \hline
        \endfoot
        \endlastfoot
&&\\
1 &
\begin{tabular}{ll} &
{\small $\gamma_{12}<0, \gamma_{13}<0, \gamma_{21}>0$,} {\small $\gamma_{23}>0, \gamma_{31}>0, \gamma_{32}<0$}\qquad \qquad \qquad \\
\end{tabular}
&
\parbox{2cm}{\vspace{2pt}\includegraphics[width=1.55cm,height=1.15cm]{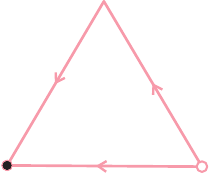}}\\
&&\\[-10pt]
    2 &
\begin{tabular}{ll} {\small  \quad \
 (i)}& {\small $\gamma_{12}<0, \gamma_{13}<0, \gamma_{21}<0$,} {\small   $\gamma_{23}>0, \gamma_{31}>0, \gamma_{32}<0$} \quad  \quad\\
{\small \quad \ (ii)}& \ {\small   $a_{31}\beta_{12}+a_{32}\beta_{21}<r_3$
}\quad \end{tabular}
 &
 \parbox{2cm}{\vspace{2pt}\includegraphics[width=1.55cm,height=1.15cm]{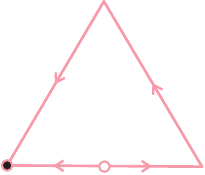}}\\
&&\\[-10pt]
    3 &
\begin{tabular}{ll} {\small
 (i)}&{\small  $\gamma_{12}<0, \gamma_{13}<0, \gamma_{21}>0$,} {\small  $\gamma_{23}<0, \gamma_{31}>0, \gamma_{32}<0$}\\
{\small  (ii)}&{\small  $a_{12}\beta_{23}+a_{13}\beta_{32}<r_1$
} \end{tabular}
 &
\parbox{2cm}{\vspace{2pt}\includegraphics[width=1.55cm,height=1.15cm]{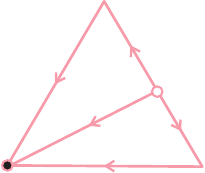}\vspace{2pt}}\\
&&\\[-10pt]
4 &
\begin{tabular}{ll} {\small
 (i)}&{\small  $\gamma_{12}>0, \gamma_{13}<0, \gamma_{21}>0$,} {\small   $\gamma_{23}<0, \gamma_{31}>0, \gamma_{32}<0$}\\
{\small  (ii)}&{\small  $a_{12}\beta_{23}+a_{13}\beta_{32}<r_1$}
{\small \quad  (iii)} \quad {\small  $a_{31}\beta_{12}+a_{32}\beta_{21}>r_3$
} \end{tabular}
 &
\parbox{2cm}{\vspace{2pt}\includegraphics[width=1.55cm,height=1.15cm]{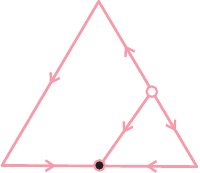}\vspace{2pt}}\\
&&\\[-10pt]
    5 &
\begin{tabular}{ll} {\small
 (i)}&{ \small $\gamma_{12}>0, \gamma_{13}>0, \gamma_{21}>0$,} {\small   $\gamma_{23}<0, \gamma_{31}<0, \gamma_{32}>0$} \qquad \\
{\small  (ii)}&{ \small $a_{31}\beta_{12}+a_{32}\beta_{21}>r_3$
} \end{tabular}
 &
    \parbox{2cm}{\vspace{2pt}\includegraphics[width=1.55cm,height=1.15cm]{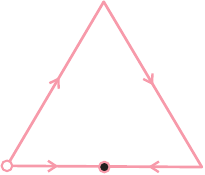}\vspace{2pt}}
\\
    &&\\[-10pt]
6 &
\begin{tabular}{ll} {\small
 (i)}&{ \small $\gamma_{12}>0, \gamma_{13}>0, \gamma_{21}<0$,} {\small   $\gamma_{23}>0, \gamma_{31}<0, \gamma_{32}>0$} \qquad \\
{\small  (ii)}&{\small \ $a_{12}\beta_{23}+a_{13}\beta_{32}>r_1$
} \end{tabular}
 &
\parbox{2cm}{\vspace{2pt}\includegraphics[width=1.55cm,height=1.15cm]{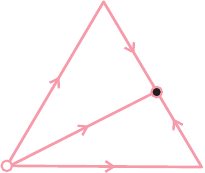}\vspace{2pt}} \\
&&\\[-10pt]
    7 &
\begin{tabular}{ll} {\small
 (i)}&{\small  $\gamma_{12}>0, \gamma_{13}>0, \gamma_{21}>0$,} { \small  $\gamma_{23}>0, \gamma_{31}<0, \gamma_{32}<0$}\\
{\small (ii)}&{\small  $a_{31}\beta_{12}+a_{32}\beta_{21}<r_3$}
 \end{tabular}
 &
\parbox{2cm}{\vspace{2pt}\includegraphics[width=1.55cm,height=1.15cm]{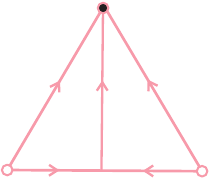}\vspace{2pt}} \\
&&\\[-10pt]
8 &
\begin{tabular}{ll} {\small
 (i)}&{\small  $\gamma_{12}>0, \gamma_{13}>0, \gamma_{21}>0$,} {\small   $\gamma_{23}<0, \gamma_{31}<0, \gamma_{32}<0$}\\
{\small  (ii)}&{\small  $a_{12}\beta_{23}+a_{13}\beta_{32}<r_1$}
{\small \quad  (iii)}\quad { \small $a_{31}\beta_{12}+a_{32}\beta_{21}<r_3$
} \end{tabular}
 &
\parbox{2cm}{\vspace{2pt}\includegraphics[width=1.55cm,height=1.15cm]{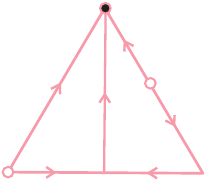}\vspace{2pt}}\\
&&\\[-10pt]
    9 &
\begin{tabular}{ll} {\small
 (i)}&{\small  $\gamma_{12}>0, \gamma_{13}>0, \gamma_{21}>0$,} {\small   $\gamma_{23}>0, \gamma_{31}<0, \gamma_{32}>0$}\\
{\small  (ii)}&{\small  $a_{12}\beta_{23}+a_{13}\beta_{32}>r_1$}
{\small \quad  (iii)}\quad {\small  $a_{31}\beta_{12}+a_{32}\beta_{21}<r_3$
} \end{tabular}
 &
\parbox{2cm}{\vspace{2pt}\includegraphics[width=1.55cm,height=1.15cm]{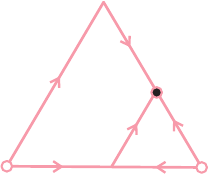}\vspace{2pt}} \\
&&\\[-10pt]
10 &
\begin{tabular}{ll} {\small
 (i)}&{\small  $\gamma_{12}>0, \gamma_{13}>0, \gamma_{21}>0$,} {\small   $\gamma_{23}>0, \gamma_{31}<0, \gamma_{32}>0$}\\
{\small  (ii)}&{\small  $a_{12}\beta_{23}+a_{13}\beta_{32}<r_1$}
{\small \quad  (iii)} \quad {\small  $a_{31}\beta_{12}+a_{32}\beta_{21}>r_3$
} \end{tabular}
 &
\parbox{2cm}{\vspace{2pt}\includegraphics[width=1.55cm,height=1.15cm]{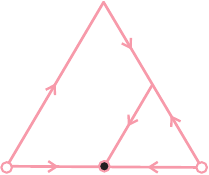}\vspace{2pt}}    \\
&&\\[-10pt]
    11 &
\begin{tabular}{ll} {\small
 (i)}&{\small  $\gamma_{12}>0, \gamma_{13}>0, \gamma_{21}>0$,} {\small   $\gamma_{23}<0, \gamma_{31}>0, \gamma_{32}<0$}\\
{\small  (ii)}& {\small  $a_{12}\beta_{23}+a_{13}\beta_{32}<r_1$}
{\small \quad  (iii)} \quad {\small  $a_{21}\beta_{13}+a_{23}\beta_{31}<r_2$}\\
{\small  (iv)}& { \small $a_{31}\beta_{12}+a_{32}\beta_{21}>r_3$
} \end{tabular}
 &
 \parbox{2cm}{\vspace{2pt}\includegraphics[width=1.55cm,height=1.15cm]{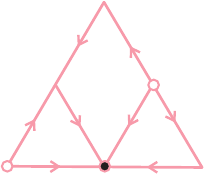}\vspace{2pt}}   \\
&&\\[-10pt]
12 &
\begin{tabular}{ll} {\small
 (i)}&{\small  $\gamma_{12}>0, \gamma_{13}>0, \gamma_{21}>0$,} {\small   $\gamma_{23}>0, \gamma_{31}>0, \gamma_{32}>0$}\\
{\small  (ii)}&{\small  $a_{12}\beta_{23}+a_{13}\beta_{32}<r_1$}
{\small \quad  (iii)}\quad {\small  $a_{21}\beta_{13}+a_{23}\beta_{31}<r_2$}\\
{\small  (iv)}&{\small  $a_{31}\beta_{12}+a_{32}\beta_{21}>r_3$
} \end{tabular}
 &
\parbox{2cm}{\vspace{2pt}\includegraphics[width=1.55cm,height=1.15cm]{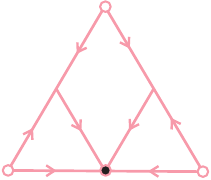}\vspace{2pt}}   \\
&&\\[-10pt]
    13 &
\begin{tabular}{ll} {\small
 (i)}&{\small $\gamma_{12}<0, \gamma_{13}<0, \gamma_{21}<0$,} {   $\gamma_{23}<0, \gamma_{31}>0, \gamma_{32}>0$}\\
{\small  (ii)}&{\small  $a_{31}\beta_{12}+a_{32}\beta_{21}>r_3$
} \end{tabular}
 &
\parbox{2cm}{\vspace{2pt}\includegraphics[width=1.55cm,height=1.15cm]{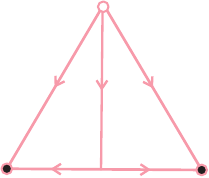}\vspace{2pt}}    \\
&&\\[-10pt]
14 &
\begin{tabular}{ll} {\small
 (i)}&{\small  $\gamma_{12}<0, \gamma_{13}<0, \gamma_{21}<0$,} {\small   $\gamma_{23}>0, \gamma_{31}>0, \gamma_{32}>0$}\\
{\small  (ii)}&{ \small $a_{12}\beta_{23}+a_{13}\beta_{32}>r_1$}
{\small \quad (iii)}\quad{ \small $a_{31}\beta_{12}+a_{32}\beta_{21}>r_3$
} \end{tabular}
 &
\parbox{2cm}{\vspace{2pt}\includegraphics[width=1.55cm,height=1.15cm]{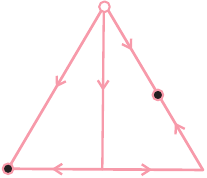}\vspace{2pt}}   \\
&&\\[-10pt]
    15 &
\begin{tabular}{ll} {\small \quad
 (i)}&{\small  $\gamma_{12}<0, \gamma_{13}<0, \gamma_{21}<0$,} {\small  $\gamma_{23}<0, \gamma_{31}>0, \gamma_{32}<0$}\\
{\small  \quad (ii)}&{$a_{12}\beta_{23}+a_{13}\beta_{32}<r_1$}
{\small \quad  (iii)}\quad{$a_{31}\beta_{12}+a_{32}\beta_{21}>r_3$
} \end{tabular}
 &
\parbox{2cm}{\vspace{2pt}\includegraphics[width=1.55cm,height=1.15cm]{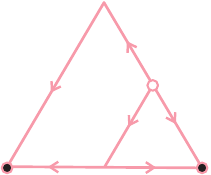}\vspace{2pt}}   \\
&&\\[-10pt]
16 &
\begin{tabular}{ll} {\small
 (i)}&{\small  $\gamma_{12}<0, \gamma_{13}<0, \gamma_{21}<0$,} {\small   $\gamma_{23}<0, \gamma_{31}>0, \gamma_{32}<0$}\\
{\small  (ii)}&{\small  $a_{12}\beta_{23}+a_{13}\beta_{32}>r_1$}
{\small \quad  (iii)}\quad {$a_{31}\beta_{12}+a_{32}\beta_{21}<r_3$
} \end{tabular}
 &
\parbox{2cm}{\vspace{2pt}\includegraphics[width=1.55cm,height=1.15cm]{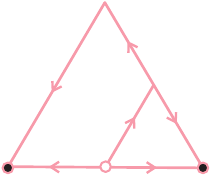}\vspace{2pt}}   \\
&&\\[-10pt]
    17 &
\begin{tabular}{ll} {\small
 (i)}&{\small  $\gamma_{12}<0, \gamma_{13}<0, \gamma_{21}<0$,} {\small   $\gamma_{23}>0, \gamma_{31}<0, \gamma_{32}>0$}\\
{\small  (ii)}&{\small  $a_{12}\beta_{23}+a_{13}\beta_{32}>r_1$}
{\small \quad  (iii)}\quad {$a_{21}\beta_{13}+a_{23}\beta_{31}>r_2$}\\
{\small  (iv)}&{\small  $a_{31}\beta_{12}+a_{32}\beta_{21}<r_3$
} \end{tabular}
 &
\parbox{2cm}{\vspace{2pt}\includegraphics[width=1.55cm,height=1.15cm]{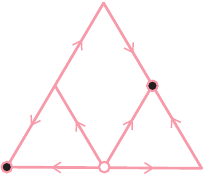}\vspace{2pt}}   \\
&&\\[-10pt]
    18 &
\begin{tabular}{ll} {\small
 (i)}&{\small  $\gamma_{12}<0, \gamma_{13}<0, \gamma_{21}<0$,} {\small   $\gamma_{23}<0, \gamma_{31}<0, \gamma_{32}<0$}\\
{\small  (ii)}&{\small  $a_{12}\beta_{23}+a_{13}\beta_{32}>r_1$}
{\small \quad  (iii)}\quad {$a_{21}\beta_{13}+a_{23}\beta_{31}>r_2$}\\
{\small  (iv)}& {\small  $a_{31}\beta_{12}+a_{32}\beta_{21}<r_3$
} \end{tabular}
 &
\parbox{2cm}{\vspace{2pt}\includegraphics[width=1.55cm,height=1.15cm]{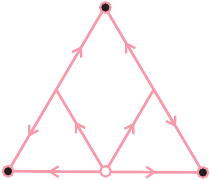}\vspace{2pt}}   \\
&&\\[-10pt]
19 &
\begin{tabular}{ll} {\small
 (i)}&{\small  $\gamma_{12}>0, \gamma_{13}>0, \gamma_{21}<0$,} {\small   $\gamma_{23}<0, \gamma_{31}<0, \gamma_{32}<0$} \qquad \\
{\small  (ii)}&{\small  $a_{12}\beta_{23}+a_{13}\beta_{32}<r_1$ \qquad
} \end{tabular}
 &
\parbox{2cm}{\vspace{2pt}\includegraphics[width=1.55cm,height=1.15cm]{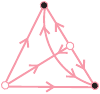}\vspace{2pt}} \\
&&\\[-10pt]
    20 &
\begin{tabular}{ll} {\small
 (i)}&{\small  $\gamma_{12}<0, \gamma_{13}<0, \gamma_{21}<0$,} {\small   $\gamma_{23}<0, \gamma_{31}>0, \gamma_{32}<0$}\\
{\small  (ii)}&{\small  $a_{12}\beta_{23}+a_{13}\beta_{32}<r_1$}
{\small \quad  (iii)}\quad{  $a_{31}\beta_{12}+a_{32}\beta_{21}<r_3$
} \end{tabular}
 &
 \parbox{2cm}{\vspace{2pt}\includegraphics[width=1.55cm,height=1.15cm]{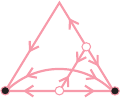}\vspace{2pt}} \\
&&\\[-10pt]
21 &
\begin{tabular}{ll} {\small
 (i)}&{\small  $\gamma_{12}<0, \gamma_{13}<0, \gamma_{21}<0$,} {\small   $\gamma_{23}>0, \gamma_{31}<0, \gamma_{32}>0$}\\
{\small  (ii)}&{\small  $a_{12}\beta_{23}+a_{13}\beta_{32}>r_1$}
{\small  \quad (iii)}\quad {  $a_{21}\beta_{13}+a_{23}\beta_{31}<r_2$}\\
{\small  (iv)}& {\small  $a_{31}\beta_{12}+a_{32}\beta_{21}<r_3$
} \end{tabular}
 &
 \parbox{2cm}{\vspace{2pt}\includegraphics[width=1.55cm,height=1.15cm]{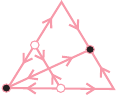}\vspace{2pt}} \\
&&\\[-10pt]
    22 &
\begin{tabular}{ll} {\small
 (i)}&{\small  $\gamma_{12}>0, \gamma_{13}>0, \gamma_{21}<0$,} {\small   $\gamma_{23}<0, \gamma_{31}>0, \gamma_{32}<0$}\\
{\small  (ii)}&{\small  $a_{12}\beta_{23}+a_{13}\beta_{32}<r_1$}
{\small \quad  (iii)}\quad { \small $a_{21}\beta_{13}+a_{23}\beta_{31}>r_2$
} \end{tabular}
 &
 \parbox{2cm}{\vspace{2pt}\includegraphics[width=1.55cm,height=1.15cm]{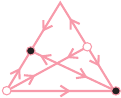}\vspace{2pt}} \\
&&\\[-10pt]
23 &
\begin{tabular}{ll} {\small
 (i)}&{\small  $\gamma_{12}>0, \gamma_{13}>0, \gamma_{21}>0$,} {\small   $\gamma_{23}>0, \gamma_{31}<0, \gamma_{32}<0$}\\
{\small  (ii)}&{\small  $a_{31}\beta_{12}+a_{32}\beta_{21}>r_3$
} \end{tabular}
 &
 \parbox{2cm}{\vspace{2pt}\includegraphics[width=1.55cm,height=1.15cm]{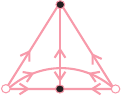}\vspace{2pt}} \\
&&\\[-10pt]
    24 &
\begin{tabular}{ll} {\small
 (i)}&{\small  $\gamma_{12}>0, \gamma_{13}>0, \gamma_{21}>0$,} {\small   $\gamma_{23}>0, \gamma_{31}<0, \gamma_{32}>0$}\\
{\small  (ii)}&{\small  $a_{12}\beta_{23}+a_{13}\beta_{32}>r_1$}
{\small \quad  (iii)}\quad {\small  $a_{31}\beta_{12}+a_{32}\beta_{21}>r_3$
} \end{tabular}
&
 \parbox{2cm}{\vspace{2pt}\includegraphics[width=1.55cm,height=1.15cm]{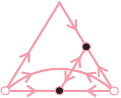}\vspace{2pt}} \\
&&\\[-10pt]
25 &
\begin{tabular}{ll} {\small
 (i)}&{\small  $\gamma_{12}>0, \gamma_{13}>0, \gamma_{21}>0$,} {\small   $\gamma_{23}<0, \gamma_{31}>0, \gamma_{32}<0$}\\
{\small  (ii)}&{\small  $a_{12}\beta_{23}+a_{13}\beta_{32}<r_1$}
{\small \quad  (iii)}\quad{\small  $a_{21}\beta_{13}+a_{23}\beta_{31}>r_2$}\\
{\small  (iv)}&{\small  $a_{31}\beta_{12}+a_{32}\beta_{21}>r_3$
} \end{tabular}
 &
 \parbox{2cm}{\vspace{2pt}\includegraphics[width=1.55cm,height=1.15cm]{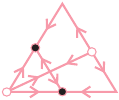}\vspace{2pt}} \\
&&\\[-10pt]
    26 (a) &
\begin{tabular}{ll} {\small
 (i)}&{\small  $\gamma_{12}>0, \gamma_{13}>0, \gamma_{21}<0$,} {\small   $\gamma_{23}<0, \gamma_{31}>0, \gamma_{32}<0$}\\
{\small  (ii)}&{\small  $a_{12}\beta_{23}+a_{13}\beta_{32}>r_1$}
{\small \quad  (iii)}\quad{ \small $a_{21}\beta_{13}+a_{23}\beta_{31}<r_2$
}\\
{\small  (iv)}&{\small  $\zeta<0$}
\end{tabular}
 &
\parbox{2cm}{\vspace{2pt}\includegraphics[width=1.55cm,height=1.15cm]{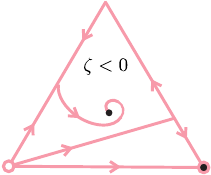}\vspace{2pt}}\\
&&\\[-10pt]
    26 (b) &
\begin{tabular}{ll} {\small
 (i)}&{\small  $\gamma_{12}>0, \gamma_{13}>0, \gamma_{21}<0$,} {\small   $\gamma_{23}<0, \gamma_{31}>0, \gamma_{32}<0$}\\
{\small  (ii)}&{\small  $a_{12}\beta_{23}+a_{13}\beta_{32}>r_1$}
{\small \quad  (iii)}\quad{\small  $a_{21}\beta_{13}+a_{23}\beta_{31}<r_2$
}\\
{\small  (iv)}&{\small  $\zeta>0$} \end{tabular}
 &
\parbox{2cm}{\vspace{2pt}\includegraphics[width=1.55cm,height=1.15cm]{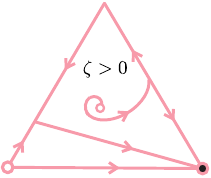}\vspace{2pt}}\\
&&\\[-10pt]
    26 (c)&
\begin{tabular}{ll} {\small
 (i)}&{\small  $\gamma_{12}>0, \gamma_{13}>0, \gamma_{21}<0$,} { \small  $\gamma_{23}<0, \gamma_{31}>0, \gamma_{32}<0$}\\
{\small  (ii)}&{\small  $a_{12}\beta_{23}+a_{13}\beta_{32}>r_1$}
{\small \quad  (iii)}\quad{\small  $a_{21}\beta_{13}+a_{23}\beta_{31}<r_2$
}\\
{\small  (iv)}&{\small  $\zeta=0$} \end{tabular}
 &
\parbox{2cm}{\vspace{2pt}\includegraphics[width=1.55cm,height=1.15cm]{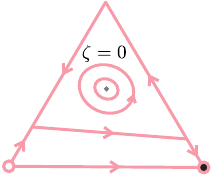}\vspace{2pt}}\\
&&\\[-10pt]
27  (a)&
\begin{tabular}{ll} {\small (i)}&{\small
 $\gamma_{12}>0, \gamma_{13}<0, \gamma_{21}<0$,} {\small   $\gamma_{23}>0, \gamma_{31}>0, \gamma_{32}<0$
}\\
{\small (ii)}& {\small $\zeta<0$}
\end{tabular}
 &
    \parbox{2cm}{\vspace{2pt}\includegraphics[width=1.55cm,height=1.15cm]{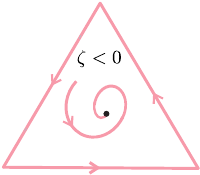}\vspace{2pt}}\\
&&\\[-10pt]
27  (b)&
\begin{tabular}{ll} {\small (i)}&{\small
 $\gamma_{12}>0, \gamma_{13}<0, \gamma_{21}<0$,} { \small  $\gamma_{23}>0, \gamma_{31}>0, \gamma_{32}<0$
}\\
{\small (ii)}& {\small $\zeta>0$}
\end{tabular}
 &
    \parbox{2cm}{\vspace{2pt}\includegraphics[width=1.55cm,height=1.15cm]{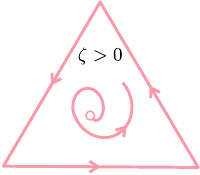}\vspace{2pt}}\\
&&\\[-10pt]
27  (c)&
\begin{tabular}{ll} {\small (i)}&{\small
 $\gamma_{12}>0, \gamma_{13}<0, \gamma_{21}<0$,} {\small   $\gamma_{23}>0, \gamma_{31}>0, \gamma_{32}<0$
}\\
{\small (ii)}& {\small $\zeta=0$}
\end{tabular}
 &
    \parbox{2cm}{\vspace{2pt}\includegraphics[width=1.55cm,height=1.15cm]{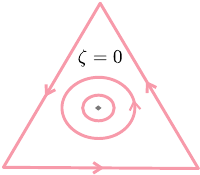}\vspace{2pt}}\\
&&\\[-10pt]
    28 &
\begin{tabular}{ll} {\small
 (i)}&{\small  $\gamma_{12}<0, \gamma_{13}<0, \gamma_{21}<0$,} {\small  $\gamma_{23}>0, \gamma_{31}>0, \gamma_{32}<0$}\\
{\small  (ii)}&{\small   $a_{31}\beta_{12}+a_{32}\beta_{21}>r_3$
} \end{tabular}
&
\parbox{2cm}{\vspace{2pt}\includegraphics[width=1.55cm,height=1.15cm]{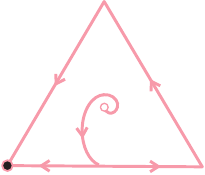}\vspace{2pt}} \\
&&\\[-10pt]
29 &
\begin{tabular}{ll} {\small
 (i)}&{\small  $\gamma_{12}>0, \gamma_{13}>0, \gamma_{21}>0$,} {\small   $\gamma_{23}<0, \gamma_{31}<0, \gamma_{32}>0$}\\
{\small  (ii)}&{\small  $a_{31}\beta_{12}+a_{32}\beta_{21}<r_3$
} \end{tabular}
 &
\parbox{2cm}{\vspace{2pt}\includegraphics[width=1.55cm,height=1.15cm]{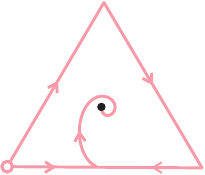}\vspace{2pt}} \\
&&\\[-10pt]
    30 &
\begin{tabular}{ll} {\small \quad \
(i)}&{\small  $\gamma_{12}<0, \gamma_{13}<0, \gamma_{21}<0$,} {\small  $\gamma_{23}<0, \gamma_{31}>0, \gamma_{32}<0$}\\
{\small \quad \  (ii)}&{\small  $a_{12}\beta_{23}+a_{13}\beta_{32}>r_1$}
{ \small \quad (iii)}\quad {  $a_{31}\beta_{12}+a_{32}\beta_{21}>r_3$
} \end{tabular}
&
 \parbox{2cm}{\vspace{2pt}\includegraphics[width=1.55cm,height=1.15cm]{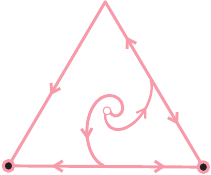}\vspace{2pt}} \\
&&\\[-10pt]
    31 &
\begin{tabular}{ll} {\small
(i)}&{\small  $\gamma_{12}>0, \gamma_{13}>0, \gamma_{21}>0$,} {\small  $\gamma_{23}>0, \gamma_{31}<0, \gamma_{32}>0$}\\
{\small  (ii)}&{\small  $a_{12}\beta_{23}+a_{13}\beta_{32}<r_1$}
{\small \quad   (iii)}\quad {\small  $a_{31}\beta_{12}+a_{32}\beta_{21}<r_3$
} \end{tabular}
&
\parbox{2cm}{\vspace{2pt}\includegraphics[width=1.55cm,height=1.15cm]{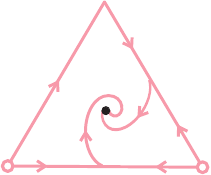}\vspace{2pt}} \\
&&\\[-10pt]
32 &
\begin{tabular}{ll} {\small
 (i)}&{\small  $\gamma_{12}<0, \gamma_{13}<0, \gamma_{21}<0$,} { \small  $\gamma_{23}<0, \gamma_{31}<0, \gamma_{32}<0$}\\
{\small  (ii)}&{\small  $a_{12}\beta_{23}+a_{13}\beta_{32}>r_1$}
{\small  \quad (iii)}\quad{\small  $a_{21}\beta_{13}+a_{23}\beta_{31}>r_2$}\\
{\small  (iv)}&{\small  $a_{31}\beta_{12}+a_{32}\beta_{21}>r_3$
} \end{tabular}
 &
 \parbox{2cm}{\vspace{2pt}\includegraphics[width=1.55cm,height=1.15cm]{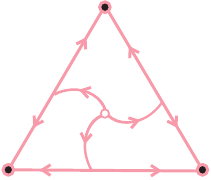}\vspace{2pt}} \\
&&\\[-10pt]
    33 &
\begin{tabular}{ll} {\small
 (i)}&{\small  $\gamma_{12}>0, \gamma_{13}>0, \gamma_{21}>0$,} {\small   $\gamma_{23}>0, \gamma_{31}>0, \gamma_{32}>0$}\\
{\small  (ii)}&{\small  $a_{12}\beta_{23}+a_{13}\beta_{32}<r_1$}
{\small \quad  (iii)}\quad {\small  $a_{21}\beta_{13}+a_{23}\beta_{31}<r_2$}\\
{\small  (iv)}&{\small  $a_{31}\beta_{12}+a_{32}\beta_{21}<r_3$
} \end{tabular}
&
 \parbox{2cm}{\vspace{2pt}\includegraphics[width=1.55cm,height=1.15cm]{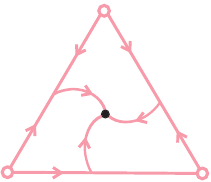}\vspace{2pt}} \\
[-8pt]\label{biao0}
\end{longtable}
\end{center}

\end{document}

%% file: ex_shared.tex
\usepackage{lipsum}
\usepackage{amsfonts}
\usepackage{graphicx}
\usepackage{epstopdf}
\usepackage{algorithmic}
\usepackage{mathtools}
\usepackage{mathrsfs}
\usepackage{xcolor}
\usepackage{enumerate}
\usepackage{amssymb}
\usepackage{subfigure}
\usepackage{longtable,booktabs}
\usepackage{lineno,hyperref}
\usepackage{threeparttable}
\newsiamthm{example}{Example}[section]

\ifpdf
  \DeclareGraphicsExtensions{.eps,.pdf,.png,.jpg}
\else
  \DeclareGraphicsExtensions{.eps}
\fi

\newsiamremark{remark}{Remark}
\newsiamremark{hypothesis}{Hypothesis}
\crefname{hypothesis}{Hypothesis}{Hypotheses}
\newsiamthm{claim}{Claim}

\headers{3-D L-V competition model with seasonal succession}{Niu, Wang and Xie}

\title{Global dynamics of three-dimensional Lotka-Volterra competition models with seasonal succession: II. Uniqueness of positive fixed points\thanks{
\textbf{Funding:} {Niu's research was supported by NSFC (grant 12001096) and the Natural Science Foundation of Shanghai (grant
24ZR1401900). Wang's research was supported by NSFC (grant 12331006), the Strategic Priority Research Program of Chinese Academy of Sciences (XDB0900100), and National Key R\&D Program of China (2024YFA1013603). Xie's research was supported by NSFC (grant 12401495) and the Natural Science Foundation of Fujian Province (grant 2026J001752).}}}

\author{Lei Niu\thanks{School of Mathematics and Statistics, Donghua University, Shanghai
201620, China (lei.niu@ \, \, dhu.edu.cn).}
\and Yi Wang\thanks{School of Mathematical Sciences, University of Science and Technology of China,
Hefei, Anhui 230026, China
  (\email{wangyi@ustc.edu.cn}).}
\and Xizhuang Xie\thanks{Corresponding author. School of Mathematical Sciences, Huaqiao University,
Quanzhou, Fujian 362021, China
(xzx@hqu.edu.cn)}.\\
\hspace{2em} \textcolor{magenta}{This manuscript has been accepted for publication in SIAM Journal on Applied Mathematics.} }

\usepackage{amsopn}
